\documentclass[12pt,epsf]{article}
\usepackage{subfigure}
\usepackage{units}
\usepackage{color}
\usepackage{ulem}
\usepackage{version}

\def\paragraph#1{{\bf #1\ }}

 \def\OO{\rm \hbox{O\kern-.34em\raise.47ex
         \hbox{$\scriptscriptstyle |$}\kern-.46em\raise.47ex
         \hbox{$\scriptscriptstyle |$}\kern+0.5 em }}
\def\RR{\mbox{\mathrm I\hspace{-0.50ex}R} }

\def\hcboxcm#1#2{\hbox to #1{\hfill #2 \hfill}}

\def\null{\hbox{}}
\let\eps\varepsilon

\def\tn1{\widetilde n_1}
\def\tn2{\widetilde n_2}
\def\tn{\widetilde n }

\let\ds\displaystyle

\def\be{\begin{equation}}
\def\ee{\end{equation}}
\def\eeq{\end{equation}}
\def\bea{\begin{eqnarray}}
\def\eea{\end{eqnarray}}
\def\bean{\begin{eqnarray*}}
\def\eean{\end{eqnarray*}}

\def\RR{{\mathrm{ I~\hspace{-1.15ex}R}}}

\def\qquad{{\quad\quad}}

\def\={\, = \, }

\def\Box{\leavevmode\vbox{\hrule
     \hbox{\vrule\kern4pt\vbox{\kern4pt}%
           \vrule}\hrule}}
\def\blackbox{\leavevmode\vrule height 5pt width 4pt depth 0pt\relax}
\catcode`@=11

\def\eqalign#1{\null\,\vcenter{\openup1\jot \m@th
   \ialign{\strut \hfil$\displaystyle{##}$ & $\displaystyle{{}##}$\hfil
      \crcr#1\crcr}}\,}
\def\eqalignrll#1{\null\,\vcenter{\openup1\jot \m@th
   \ialign{\strut \hfil$\displaystyle{##}$ & $\displaystyle{{}##}$\hfil
    & $\displaystyle{{}##}$\hfil
      \crcr#1\crcr}}\,}
\def\eqalignrcl#1{\null\,\vcenter{\openup1\jot \m@th
   \ialign{\strut \hfil$\displaystyle{##}$ &\hfil $\displaystyle{{}##}$\hfil
    & $\displaystyle{{}##}$\hfil
      \crcr#1\crcr}}\,}
\def\eqalignno#1{\displ@y \tabskip\@centering
  \halign to\displaywidth{\hfil$\@lign\displaystyle{##}$\tabskip\z@skip
    &$\@lign\displaystyle{{}##}$\hfil\tabskip\@centering
    &\llap{$\@lign##$}\tabskip\z@skip\crcr
    #1\crcr}}
\newcounter{appendix}
\newcounter{sectionz}
\def\appendix{\advance\c@appendix by 1
   \def\thesectionz {\Alph{appendix}}
\def\thesection{\Alph{appendix}} 
   \ifnum\c@appendix=1 \setcounter{section}{-1} \fi
   \@startsection {section}{1}{\z@}{-3.5ex plus -1ex minus 
  -.2ex}{2.3ex plus .2ex}{\large\bf}}

\catcode`@=12
\newtheorem{lemme}{Lemma}[section]  

\newtheorem{theorem}[lemme]{Theorem}

\newtheorem{corollary}[lemme]{Corollary}

\newtheorem{definition}[lemme]{Definition}

\newtheorem{proposition}[lemme]{Proposition}

\newtheorem{remark}[lemme]{Remark} 

\def\deblem{\begin{lemme}\it }
\def\finlem{\end{lemme}}
\def\debthm{\begin{theorem}\it }
\def\finthm{\end{theorem}}
\def\debprop{\begin{proposition} \it}
\def\finprop{\end{proposition}}
\def\debcor{\begin{corollary}\it }
\def\fincor{\end{corollary}}
\def\debdef{\begin{definition}\it}
\def\findef{\end{definition}}
\def\debrem{\begin{remark}\em}
\def\finrem{\null\hfill\blackbox\end{remark}}
\def\debproof{{\bf Proof: \ }}
\def\finproof{\null\hfill {$\blackbox$}\bigskip}  
\def\OO{\mathbb{O}}

\def\RR{\mathbb{R}}

\usepackage{amsfonts}
\usepackage{amssymb}
\usepackage{amsbsy}
\usepackage{amsmath}

\usepackage{stmaryrd}
\usepackage[dvips]{graphicx}
\usepackage{psfrag}
\usepackage[hang,center]{caption}
\usepackage{float}
\usepackage{cite}

\usepackage{latexsym, amssymb, enumerate, amsmath}
\usepackage{bbm}
\usepackage{epsfig}
\usepackage{subfigure}

\usepackage{multirow}
\usepackage[affil-it]{authblk}

\begin{document}

\title{A hybrid method for anisotropic elliptic problems based on the coupling of an Asymptotic-Preserving method with the Asymptotic-Limit model.}

\author{A. Crestetto, F. Deluzet, C. Negulescu\thanks{Electronic address: \texttt{anais.crestetto@math.univ-toulouse.fr; fabrice.deluzet@math.univ-toulouse.fr; claudia.negulescu@math.univ-toulouse.fr}}}

\affil{Institut de Math\'ematiques de Toulouse; UMR 5219\\
Universit\'e de Toulouse; CNRS\\
UPS IMT, 118 route de Narbonne; F-31062 Toulouse, France}

%%%%%%%%%%%%%%%%%%%%%%%%%%%%%%%%%%%%%%%%%%%%%%%%%%%%%%%%%%%%%%%%%%%%%%%%%%%%%%%
%\begin{document}

\maketitle

\renewcommand{\thefootnote}{\fnsymbol{footnote}}

\renewcommand{\thefootnote}{\arabic{footnote}}

\begin{abstract}
This paper presents a hybrid numerical method to solve efficiently a class of highly anisotropic elliptic problems. The anisotropy is aligned with one coordinate-axis and its strength is described by a parameter $\eps \in (0,1]$, which can largely vary in the study domain. Our hybrid model is based on asymptotic techniques and couples (spatially) an Asymptotic-Preserving model with its asymptotic Limit model, the latter being used in regions where the anisotropy parameter $\eps$ is small.  Adequate coupling conditions link the two models. Aim of this hybrid procedure is to reduce the computational time for problems where the region of small $\eps$-values extends over a significant part of the domain, and this due to the reduced complexity of the limit model.

\end{abstract}

\bigskip


%\tableofcontents

%%%%%%%%%%%%%%%%%%%%%%%%%%%%%%%%%%%%%%%%%%%%%%
\section{Introduction} \label{SEC1}
%%%%%%%%%%%%%%%%%%%%%%%%%%%%%%%%%%%%%%%%%%%%%%%

The present work is a contribution to the numerical resolution of highly anisotropic elliptic equations, the anisotropy being aligned with one coordinate axis and described by a perturbation parameter $\eps \in (0,1]$, varying considerably in the study domain. The approach presented here is based on a coupling strategy, solving an Asymptotic-Preserving reformulation of the elliptic problem where $\eps$ is non-negligible and solving the corresponding asymptotic Limit model, where $\eps$ is quasi vanishing. The  strategy we propose is particularly well suited for physical systems in which the anisotropy parameter $\varepsilon$ is very small in a large part of the study domain.

Such kind of directionally anisotropic diffusion systems are common in physical applications, such as plasma physics \cite{chen,GoldsRuth,heikkila,Jardin,schunk,stix}. 
%Such kind of directionally anisotropic diffusion systems are common in many physical applications, such as medical diagnostics via Diffusion Tensor Magnetic Resonance Imaging \cite{basser}, flows in porous media \cite{Bear},  image processing \cite{weickert}, semiconductor modeling \cite{manku} and plasma physics \cite{stix}. 

The application which was at the origin of the present work comes from strongly magnetized ionospheric plasmas \cite{BCDDGT,Sami3, Sami3-2}. The problem we shall study here is extracted from the Dynamo model and represents an elliptic equation for the computation of the electric potential in 2D, {\it i.e.}
\begin{equation}
\label{eq:SP}
(P)\,\,\, \left\{\begin{array}{ll}
 -\nabla\cdot\left(\mathbb{A}\nabla u\right)=f, & \text{~in~} \Omega:=\Omega_x\times\Omega_z, \\
(\mathbb{A}\, \nabla u) \cdot n = g, & \text{~on~} \Omega_x \times \partial\Omega_z, \\
 u=0, & \text{~on~} \partial\Omega_x \times \Omega_z,
\end{array}
\right.
\end{equation}
where $\Omega_x\subset\mathbb{R}$, $\Omega_z\subset\mathbb{R}$ are intervals and $\partial \Omega$ denotes the boundary of $\Omega$, with outward normal $n$. We assume that the anisotropy direction is fixed and aligned with the $z$-coordinate, the diffusion matrix $\mathbb{A}$ being thus given by
\begin{equation}\label{eq:SP:tensor}
\mathbb{A}=\left(\begin{array}{cc}
A_x & 0\\
0 & \frac{1}{\varepsilon}A_z
\end{array}\right),
\end{equation}
with $A_x$ and $A_z$ of the same order of magnitude. The high anisotropy of the problem is parametrized by $\varepsilon\in (0,1]$ that can become very small in some regions of $\Omega$. In the ionospheric plasma framework referred above, this parameter $\varepsilon$ is the ratio of the collision frequency to the cyclotron frequency. The aim of this paper is to propose a domain decomposition strategy for an efficient numerical resolution of this singularly perturbed system (\ref{eq:SP})-(\ref{eq:SP:tensor}) (P-model) in situations where $\varepsilon$ undergoes large variations along the $z$ coordinate. Generalizations of this strategy to 3D problems, where ${\mathbf x}=(x,y) \in \Omega_{\mathbf x} \subset \RR^2$, $A_{\mathbf x}({\mathbf x},z) \in \RR^{ 2 \times 2}$, while the anisotropy is  remaining aligned with the $z$-direction, are straightforward.\\

A naive resolution of (\ref{eq:SP})-(\ref{eq:SP:tensor}) leads to an unusable scheme in the limit $\varepsilon\rightarrow 0$, because (\ref{eq:SP})-(\ref{eq:SP:tensor}) degenerates into an ill-posed problem. Indeed, supposing $\eps$ constant and letting it tend formally towards zero, yields the reduced model
\be \label{R}
\left(R\right)\,\,\,
\left\{
\begin{array}{ll}
  -\partial_z \left( A_z  \partial_z u \right)=0, & \textrm{ for } \left(x,z\right) \in \Omega_x \times \Omega_z,\\[1mm]
 \partial_z u=0, &  \text{~on~} \Omega_x \times \partial\Omega_z,\\[1mm]
 u =0, & \text{~on~} \partial\Omega_x \times \Omega_z.
\end{array}
\right.
\ee
All functions which are constant along the $z$-coordinate and which satisfy the boundary condition on $\partial\Omega_x \times \Omega_z$ are solutions of this ill-posed R-model. This non-uniqueness leads to an ill-conditioned linear system for the discretized P-model, when $\varepsilon\rightarrow 0$. To avoid this degeneracy, one way is to use an asymptotic-preserving AP-scheme, notion introduced firstly by Jin in \cite{APJin}. These schemes are based on asymptotic techniques and consist in a reformulation of the singularly-perturbed problem $(P)$ into an equivalent problem $(AP)$, which in the limit $\eps \rightarrow 0$ yields the ``right'' Limit-model, defined as the problem satisfied by the $\eps \rightarrow 0$ limit of the solution-sequence $\{u^\eps\}_{\eps>0}$. In this case, the L-model is given by
$$
\left(L\right) 
\left\{
\begin{array}{ll}
 \ds - \partial_x \left( \overline{A_x}  \partial_x {u}_0\right) =\overline{f} + { g_+ \over L_z} -{ g_- \over  L_z }, & \textrm{ for }  x \in \Omega_x,\\[2mm]
 \ds  {u}_0\left(x_\pm\right)=0\,,
\end{array}
\right.
$$
where the bars signify the average over the $z$-direction (anisotropy-direction). We refer the reader to Section \ref{SSEC22}  for its derivation. This AP-procedure permits to solve special types of singularly-perturbed problems, equally accurate with respect to $\eps$, and in particular even its limit case $\eps =0$. In recent works, such AP-schemes have been developed for highly anisotropic elliptic problems of type \eqref{eq:SP}. In \cite{DDN}, Degond \textit{et al.} present an AP-scheme for the $\varepsilon$-constant case, based on the decomposition of the unknown $u$ into its mean part along the $z$-direction and the fluctuation part (so-called duality-based strategy). The reformulated system is solved iteratively, starting from an approximation of the fluctuation. A direct resolution is proposed in \cite{BDNY,PhDYang}, as well as a generalization to a variable anisotropy strength $\varepsilon\left(z\right)$, presenting steep gradients. A further extension to a variable anisotropy direction, given by a known vector field, is treated in  \cite{DDLNN}, using in addition Lagrange multiplier techniques.\\ 

Aim of this paper is to increase the efficiency of the duality-based method of \cite{BDNY} (anisotropy aligned with one coordinate axis), in the particular case of an anisotropy-strength $\eps(z)$ which is very small in a large part of the domain. A full AP-scheme like in \cite{BDNY,DDN} would be one possibility to solve such kind of problem. It appears however that in the region of small $\eps$-values, solving the Limit-model is computationally more interesting in the framework of an anisotropy direction aligned with one coordinate axis. The reason for this is that the Limit-model is a lower-dimensional problem, in the present case a 1D, $z$-independent, elliptic problem.

These considerations lead us to the introduction of a domain decomposition strategy, where the AP-model is used where $\varepsilon\left(z\right)$ is of order one, and the L-model in the regions where $\varepsilon\left(z\right)$ is "small" enough, both models being coupled with appropriate interface conditions. Such a coupling is studied in the one-dimensional framework in \cite{DDMNNP}. We propose here to extend this coupling strategy to 2D problems and to provide the first analysis results. \\

Domain decomposition methods \cite{QUA} are standard techniques to derive computationally efficient numerical schemes for problems showing different behaviors in different regions of the domain.  For example in many-particle dynamics, microscopic models (Boltzmann equations, corresponding to our P-model) are coupled to macroscopic models (Euler equations, corresponding to our L-model) via appropriate interface conditions \cite{klar1,klar2,Sch}. Usually the P-model is singularly perturbed, meaning that its asymptotic limit model is of different nature (e.g. kinetic compared to fluid models), such that the coupling of both models is very delicate. Indeed, the positioning of the interface requires the existence of a zone, where both models are valid and an automatic detection criteria. Moreover  the design and computation of appropriate coupling conditions  is generally  challenging from an analytic as well as numerical point of view. In the transition from a kinetic equation to its hydrodynamic or diffusion limit, an issue has been proposed in \cite{DJ, Luc} to avoid interface conditions. It consists in using a buffer zone, where both models have to be solved. However, it is not always possible to find such a buffer zone, where both models are accurate. Our (AP/L)-coupling strategy overcomes all these difficulties, thanks to the fact that the L-model is automatically recovered from the AP formulation as $\varepsilon\rightarrow 0$. In other words, the AP-model is a regular perturbation of the L-model, fact which permits their coupling via simple interface conditions. Moreover, the AP-formulation gives accurate results for all values $\varepsilon\in (0,1]$, such that the positioning of the interface is no more a problem.\\

~

The outline of this paper is the following. Section~\ref{SEC2} presents the Singular-Pertur\-bation P-problem, its asymptotic limit L-model and an Asymp\-totic-Preserving AP-reformulation, which was introduced in previous works (see \cite{BDNY,DDLNN,DDN}). Section~\ref{SEC3} is devoted to the coupling of the AP-reformulation with the L-model via Dirichlet-Neumann transfer conditions. We first explain this coupling-strategy and then analyze it rigorously. We present in Section~\ref{SEC4} the numerical discretization based on a finite element method and comment the obtained numerical results. Section~\ref{SEC5} is devoted to some conclusions.% and perspectives.

%%%%%%%%%%%%%%%%%%%%%%%%%%%%%%%%%%%%%%%%%%%%%%
\section{The elliptic problem and its Asymptotic-Preser\-ving reformulation}\label{SEC2}
%%%%%%%%%%%%%%%%%%%%%%%%%%%%%%%%%%%%%%%%%%%%%%

The anisotropic, two dimensional elliptic problem we shall consider, is posed for simplicity reasons on a rectangular domain $\Omega=\Omega_x\times\Omega_z$, where $\Omega_x:=(x_-,x_+)\subset\mathbb{R}$ and $\Omega_z:=(z_-,z_+)\subset\mathbb{R}$ are two intervals of respective length $L_x$ and $L_z$. 

Aim of this section is firstly to present the Singular-Perturbation problem we are interested in,  then to recall the corresponding asymptotic limit model as well as an Asymptotic-Preserving reformulation, introduced in some previous works  \cite{BDNY,DDLNN}. For more mathematical details, we refer the reader to these works.

%%%%%%%%%%%%%%%%%%%%%%%%%%%%%%%%%%%%%%%%%%%%%%%
\subsection{The Singular-Perturbation problem: P-model} \label{SSEC21}
%%%%%%%%%%%%%%%%%%%%%%%%%%%%%%%%%%%%%%%%%%%%%%%
We are interested in an efficient resolution of the following Singular-Perturbation problem
\be \label{P}
\left(P\right)
\left\{
\begin{array}{ll}
 - \partial_x \left( A_x \partial_x u_{\varepsilon} \right) -\partial_z \left( {A_z \over
    \eps\left(z\right)} \partial_z u_{\varepsilon} \right)
=f, & \textrm{ for }  \left(x,z\right) \in \Omega_x \times \Omega_z,\\[3mm]
  {A_z\left(x,z_{\pm}\right) \over \eps\left( z_{\pm}\right)}  \partial_z u_{\varepsilon}\left(x, z_{\pm}\right)=
g_{\pm}\left(x\right), & \textrm{ for }  x \in \Omega_x,\\[3mm]
 u_{\varepsilon}\left(x_\pm,z\right)=0, & \textrm{ for } z\in \Omega_z,
\end{array}
\right.
\ee
where we suppose that the coefficients and source terms satisfy the following hypotheses.

\vspace{0.2cm}

{\bf Hypothesis A:} {\it We consider $A_x,A_z\in L^{\infty}\left(\Omega\right)$, $f\in L^2\left(\Omega\right)$, $g_\pm\in L^2\left(\Omega_x\right)$ and $\eps \in L^\infty\left(\Omega_z\right)$, satisfying
$$
0< m_x\le A_x\left(x,z\right) \le M_x\,, \quad 0< m_z\le A_z\left(x,z\right) \le M_z\,, \quad 0 < \eps_{min} \le \eps\left(z\right) \le \eps_{max}  \le 1\,,
$$
with $m_x,m_z,M_x,M_z,\eps_{min},\eps_{max}$ some given positive constants.
}

\vspace{0.2cm}

For the mathematical investigations of the coupling strategy presented in this paper, we shall assume more regularity on the diffusion matrices and source terms, {\it i.e.}\\

\vspace{0.2cm}

{\bf Hypothesis B:} {\it Additionally to Hypothesis A we shall assume that $A_x\,, A_z\in W^{1,\infty}\left(\Omega\right)$, $(f,g_\pm) \in H^1(\Omega) \times  H^1(\Omega_x)$ and $\eps\in W^{1,\infty}\left(\Omega_z\right)$ is strictly increasing, with the bound $||\eps||_{W^{1,\infty}(\Omega_z)} \le \eps_M$.}

\vspace{0.2cm}

Let us now introduce the Hilbert-space
$$
\mathcal{V}:=\{ \psi \in H^1\left(\Omega\right) ~/~  \psi\left(x_\pm ,z\right)=0 \text{~for~} z\in \Omega_z\},
$$
associated with the scalar product
$$
\left(\phi,\psi\right)_{\mathcal{V}}:=\left(\partial_x\phi,\partial_x\psi\right)_{L^2\left(\Omega\right)}+\left(\partial_z\phi,\partial_z\psi\right)_{L^2\left(\Omega\right)}.
$$
To simplify the notations, we shall denote in the following by $(\cdot,\cdot)_{L^2}$ resp. $(\cdot,\cdot)_{L^2_x}$ the corresponding scalar-products in $L^2(\Omega)$ resp. $L^2(\Omega_x)$. 
The variational formulation of problem (P) reads then: 

\textit{Find $u_{\varepsilon} \in
\mathcal{V}$, such that}
\be \label{P_var}
\begin{array}{l}
\ds  \int_\Omega A_x\, \partial_x u_{\varepsilon} \, \partial_x \psi\,\,  \mathrm{d}x\mathrm{d}z +\int_\Omega {A_z \over
    \eps\left(z\right)} \, \partial_z u_{\varepsilon}\,  \partial_z \psi\,\, \mathrm{d}x\mathrm{d}z     \\[3mm]
\ds  ~~~~~=\left(f,\psi\right)_{L^2}+\left(g_+,\psi\left(\cdot,z_+\right)\right)_{L^2_x}-\left(g_-,\psi\left(\cdot,z_-\right)\right)_{L^2_x}~~~
    \forall \psi \in \mathcal{V}. 
\end{array}
\ee

Under Hypothesis A, this model admits a unique solution (Lax-Milgram), however its numerical approximation can not be computed for $\eps \ll 1$, by simply discretizing (\ref{P}). Indeed as explained in the introduction, (P) degenerates for $\varepsilon\rightarrow 0$ into an ill-posed problem (R), called reduced model, leading to the inversion of an ill-conditioned linear system. 
This degeneracy motivates the construction of an Asymptotic-Preserving reformulation of (P), which shall permit to get automatically the "right" asymptotic limit, when $\eps$ is vanishing. To do this, we shall first identify this limit model.

%%%%%%%%%%%%%%%%%%%%%%%%%%%%%%%%%%%%%%%%%%%%%%%
\subsection{The limit model: L-model} \label{SSEC22}
%%%%%%%%%%%%%%%%%%%%%%%%%%%%%%%%%%%%%%%%%%%%%%%

Let us consider in this subsection that $\varepsilon$ is constant. 
In the limit $\varepsilon\rightarrow 0$, the solution of the singular-perturbation model (P) is shown to converge towards some function $u_0$, solution of a limit model to be identified here. This solution is a particular solution of the R-model and it is thus independent of the $z$ variable. 

To simplify the following computations, let us introduce the following notations: for a
function $f$, we denote  by $\overline{f}$ the average over the
$z$-direction (anisotropy direction) and by $f^\prime$ the fluctuation part, given by
$$
\overline{f}\left(x\right):= {1 \over L_z} \int_{\Omega_z} f\left(x,z\right) \mathrm{d}z,~~~ f^\prime:=f-\overline{f}.
$$
We then have the following properties:
\be
\begin{array}{lllll}
\overline{f^\prime} = 0,&~~~~~&  \overline{ \left( \frac{\partial f}{\partial x} \right) }= \frac{\partial \overline{f}}{\partial x},
 &~~~~~    & 
 \overline{fg} = \overline{f} \overline{g} + \overline{f^\prime g^\prime},\\[2mm]
 \frac{\partial f}{\partial z} = \frac{\partial f^\prime}{\partial z}, &~~~~~& \left( \frac{\partial f}{\partial x} \right)^\prime = \partial \frac{f^\prime}{\partial x},&~~~~~  &\  \left(fg\right)^\prime = f^\prime g^\prime - \overline{f^\prime g^\prime}+ \overline{f} g^\prime + f^\prime \overline{g}. 
\end{array}
\ee

Now, integrating (\ref{P}) along the $z$-coordinate, passing to the limit $\eps \rightarrow 0$ and assuming that $u_\eps \rightarrow_{\eps \rightarrow 0} u_0$ in $H^1(\Omega)$ with $u_0=u_0\left(x\right)$, permits to obtain the system satisfied by the limit solution $u_0$, {\it i.e.}
\be \label{L}
\left(L\right) 
\left\{
\begin{array}{ll}
 \ds - \partial_x \left( \overline{A_x}  \partial_x {u}_0\right) =\overline{f} + { g_+ \over L_z} -{ g_- \over  L_z }, & \textrm{ for }  x \in \Omega_x,\\[2mm]
 \ds {u}_0\left(x_\pm\right)=0.
\end{array}
\right.
\ee
This so-called L-model is well-posed due to the Lax-Milgram theorem and provides an accurate solution of the P-model for very small values of $\varepsilon$. \\
The basic idea of Asymptotic-Preserving schemes  is now to reformulate the P-model in such a manner to lead automatically towards the L-model in the limit $\varepsilon\rightarrow 0$, and not towards the ill-posed R-model (\ref{R}). This procedure seems reasonable if one wants to treat, with no huge computational costs, problems with highly variable anisotropies $\eps$ within the domain.

%%%%%%%%%%%%%%%%%%%%%%%%%%%%%%%%%%%%%%%%%%%%%%%
\subsection{The fully Asymptotic-Preserving reformulation: AP-model} \label{SSEC23}
%%%%%%%%%%%%%%%%%%%%%%%%%%%%%%%%%%%%%%%%%%%%%%%
Decomposing each quantity of the P-model in its average and fluctuation part, $u_{\varepsilon}\left(x,z\right) =
\overline{u}_{\varepsilon}\left(x\right) + u_{\varepsilon}^\prime\left(x,z\right)$, permits to get an equivalent reformulation
of the Singularly-Perturbed problem (P), called Asymptotic-Preserving
model (AP-model). This model consists of two sets of equations, one for each part of the solution. Its derivation is recalled in the next lines. \\
Taking first the mean of the P-model
over $z$, we obtain the following problem, to be solved for $\overline{u}_{\varepsilon}\left(x\right)$, problem which depends on $\eps$ only through the source term 
\be \label{AP_bar}
\left(\overline{AP}\right) 
\left\{
\begin{array}{ll}
\ds   - \partial_x \left( \overline{A_x}  \partial_x \overline{u}_\eps\right) 
=\overline{f} + { g_+\over  L_z} -{ g_- \over L_z } +\partial_x \left( \overline{ A_x^\prime \partial_x u^\prime_\eps }\right),   & \textrm{ for } x \in \Omega_x,\\[2mm]
\ds \overline{u}_\eps\left(x_\pm\right)=0.
\end{array}
\right.
\ee
Secondly, decomposing now simply $u_{\varepsilon}\left(x,z\right) =
\overline{u}_{\varepsilon}\left(x\right) + u_{\varepsilon}^\prime\left(x,z\right)$ in (P), yields the following equation for $u_{\varepsilon}^\prime\left(x,z\right)$:
\be \label{AP_fluc}
\left(AP^\prime\right) 
\left\{
\begin{array}{ll}
\ds - \partial_x \left( A_x  \partial_x u^\prime_\eps \right) -\partial_z \left( {A_z \over
    \eps\left(z\right)}  \partial_z u^\prime_\eps \right)
=f + \partial_x \left( A_x  \partial_x \overline{u}_\eps
\right) , & \textrm{ for } \left(x,z\right) \in \Omega_x \times \Omega_z\\[2mm]
\ds  {A_z\left(x,z_{\pm}\right) \over \eps\left( z_{\pm}\right)} \partial_z u^\prime _\eps\left(x, z_{\pm}\right)=
g_{\pm}\left(x\right), & \textrm{ for }  x \in \Omega_x,\\[2mm]
\ds  u^\prime_\eps\left(x_\pm,z\right)=0, & \textrm{ for }  z\in  \Omega_z,\\[2mm]
 \overline{u^\prime}_\eps=0,& \textrm{ for } x\in\Omega_x \textrm{ (constraint)}.
\end{array}
\right.
\ee
The coupled system $\left(\overline{AP}\right)-\left(AP^\prime\right)$ is the Duality-Based AP-reformulation introduced
in \cite{DDLNN} and shown to be completely equivalent to the P-system for fixed $\eps>0$, leading however in the limit $\eps \rightarrow 0$ towards the well-posed L-model. The important ingredient in this reformulation is the constraint $\overline{u^\prime}_\eps=0$, which is automatically satisfied for $\varepsilon>0$, but helps to get in the limit $\varepsilon\rightarrow 0$ the constraint $\overline{u^\prime}_0=0$, which is the missing information in (R) to get the well-posed L-model. For more details we refer the reader to \cite{DDLNN}, in particular for the rigorous existence, uniqueness and $\varepsilon\rightarrow 0$ convergence proofs.\\

Let us introduce the following Hilbert-space
\begin{equation*}
\mathcal{W}:=\left\{\psi\in H^1\left(\Omega_x\right)~/~\psi\left(x_\pm \right)=0 \right\}\,,\quad 
\left(\phi,\psi\right)_\mathcal{W}:=\left(\partial_x\phi,\partial_x\psi\right)_{L^2}.
\end{equation*}
Introducing a Lagrangian multiplier  in order to cope with the constraint $\overline{u^\prime}_\eps=0$,  the variational formulation of problem $\left(\overline{AP}\right)-\left(AP^\prime\right)$ writes:

\textit{Find ($\overline{u}_{\varepsilon},u_{\varepsilon}^\prime,\overline{P}) \in \mathcal{W} \times \mathcal{V} \times L^2\left(\Omega_x\right)$ such that}
\begin{equation}
\left(AP\right)\left\{
\begin{array}{l}
\ds  (\overline{A_x}\, \partial_x\overline{u}_{\varepsilon}, \partial_x\overline{\psi})_{L^2_x}=(\overline{f},\overline{\psi})_{L^2_x}+\frac{1}{L_z} (g_+-g_-,\overline{\psi})_{L^2_x}-\frac{1}{L_z}(A_x^\prime\, \partial_xu_{\varepsilon}^\prime,\partial_x\overline{\psi})_{L^2}\\[2mm]
\ds  \hspace{12.0cm} \forall \overline{\psi}\in\mathcal{W},\\%[2mm]
\ds  (A_x\,\partial_x u_{\varepsilon}^\prime,\partial_x\psi^\prime)_{L^2}+(\frac{A_z}{\varepsilon}\,\partial_zu_{\varepsilon}^\prime,\partial_z\psi^\prime)_{L^2}+L_z\, (\overline{P}, \overline{\frac{1}{\varepsilon}\, \psi^\prime})_{L^2_x}=(f,\psi^\prime)_{L^2}\\[2mm]
\ds  \hspace{1.3cm}+(g_+,\psi^\prime(\cdot,z_+))_{L^2_x}-(g_-,\psi^\prime(\cdot,z_-))_{L^2_x}-(A_x\,\partial_x\overline{u}_{\varepsilon},\partial_x\psi^\prime)_{L^2}, \qquad \forall \psi^\prime\in\mathcal{V},\\[4mm]
\ds (\overline{Q}, \overline{u_{\varepsilon}^\prime})_{L^2_x}=0, \qquad \forall \overline{Q}\in L^2\left(\Omega_x\right).
\end{array}
\right.
\label{weakformAP}
\end{equation}

Solving (\ref{weakformAP}) instead of (\ref{P}) is numerically much more appropriate, as one gets automatically the limit problem (\ref{L}) for vanishing $\eps$, which means (\ref{weakformAP}) is not degenerating for $\eps \rightarrow 0$.
However, if the computational domain contains a large region, where the value of the
parameter $\eps$ is very small, it can be more efficient from
a computational point of view (simulation time and memory storage) to solve there directly the L-model, the latter model being of lower dimension. These considerations lead naturally to the idea of a coupling strategy between the AP-model (\ref{weakformAP}) (in regions where $\varepsilon\sim\mathcal{O}\left(1\right)$) and the 1D L-model (\ref{L}) (where $\varepsilon\ll 1$). 
%%%%%%%%%%%%%%%%%%%%%%%%%%%%%%%%%%%%%%%%%%%%%%%
\section{The (AP/L)-coupling} \label{SEC3}
%%%%%%%%%%%%%%%%%%%%%%%%%%%%%%%%%%%%%%%%%%%%%%%

This section is devoted to the coupling of the Asymptotic-Preserving reformulation (\ref{AP_bar})-(\ref{AP_fluc}) with the asymptotic limit model (\ref{L}). Firstly we shall explain the strategy and then analyze it rigorously. For simplicity reasons, we shall omit in the following the index $\varepsilon$ of the unknown $u_{\varepsilon}$.\\

Let us remark here that a coupling of the P-model with the L-model, which is usually done in literature \cite{klar1, klar2, Sch}, is not always possible. Indeed, one important issue in developing a method implementing a standard discretization of the (P)-model in a sub-domain and the (L)-model elsewhere requires that both models are solved in a common region containing the coupling interface. However, the limit model defines an accurate approximation for small values of the asymptotic parameter while the (P)-model suffers from significant precision discrepancies for these same parameter values. Defining a region where both models can produce an accurate approximation of the solution at the same time may thus be impossible (see \cite[Figure~3.3]{DDLNN} for an illustration). %These two requirements are generally met with great difficulty, and there may be even configurations with regions where none of these two problems ((P) and (L)) can be used (see discussion of section~\ref{sec:interface:choice}). 
In contrast to this, the AP-formulation gives accurate results for all values of $\eps \in [0,1]$, such that a coupling is always possible. Furthermore, due to the fact that the AP-formulation is a regular perturbation of the L-model, coupling interface conditions are easy to design. All these arguments underline the importance of the here presented strategy.
%%%%%%%%%%%%%%%%%%%
\subsection{Presentation of the coupling strategy}\label{SSEC31}
\begin{figure}[H]
\begin{center}
\psfrag{Emax}[r][r][2.5]{$\varepsilon_{min}$}
\psfrag{Emin}[r][r][2.5]{$\varepsilon_{max}$}
\psfrag{Y1}[r][r][2.]{$0$}
\psfrag{Y0}[r][r][2.]{$1$}
\psfrag{Zm1}[][][2.5]{$z_-$}
\psfrag{Zp1}[][][2.5]{$z_+$}
\includegraphics[angle=0,width=0.5\textwidth]{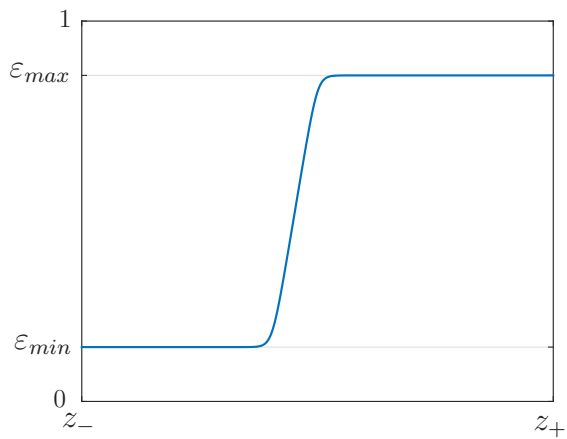}
\caption{Illustration of the anisotropy ratio $\varepsilon$ as a function of $z$.}
\label{fig:epsilonschema}
\end{center}
\end{figure}

%%%%%%%%%%%%%%%%%%%
Let us assume that in a large region of the computational domain the anisotropy parameter $\varepsilon$ is very small. For simplicity reasons, we assume that $\eps$ is a strictly increasing function and that the domain is decomposed in the $z$-direction into two sub-domains, delimited by an interface $z_\iota\in\left(z_-,z_+\right)$. Let us introduce the  following notation $\Omega_z^1:=\left(z_\iota,z_+\right)$, $\Omega_z^2:=\left(z_-,z_\iota\right)$,  $\Omega_1:= \Omega_x
\times \Omega_z^1$ and $\Omega_2:= \Omega_x \times \Omega_z^2$, and  assume $\varepsilon$ to be small in $\Omega_2$, as illustrated in Figure~\ref{fig:epsilonschema}.  Recall however that $\eps(\cdot) \in W^{1,\infty}(\Omega_z)$, such that the derivatives of $\eps(\cdot)$ have to be bounded, even when $\eps(\cdot)$ goes to zero in some parts of the domain. Our approach  is now to use in $\Omega_2$ the Limit-model and elsewhere the AP-model, with the aim to reduce the overall computational costs for the resolution of \eqref{eq:SP}. Let us detail here the coupling strategy and its variational formulation.

\subsubsection{The coupling strategy}
Let us decompose the unknown $u$ as follows $u=\overline{u}+ u^\prime_1+u^\prime_2$,
where 
$$ u^\prime_1:= u^\prime \mathbbm{1}_{\Omega_x \times \Omega_z^1}~~~\textrm{and}~~~ u^\prime_2:= u^\prime \mathbbm{1}_{\Omega_x \times \Omega_z^2}$$
are the restrictions of the fluctuation on the two sub-domains.
As mentioned above, we shall suppose that in $\Omega_2$ the parameter $\eps$ is very small, meaning that $\eps\left(z_{\iota}\right) \ll 1$. Let us now rewrite the AP-model for this new decomposition $\left(\overline{u}, u^\prime_1,u^\prime_2\right)$ coupling the two sub-domains via Dirichlet-Neumann transfer conditions, which translate the fact that the solution as well as its normal derivative $\partial_zu$ are continuous at the interface $z_{\iota}$. This new system, called in the following (AP/AP)-reformulation is completely equivalent to the AP-model ({\ref{AP_bar})-({\ref{AP_fluc}) and is given by

\begin{equation}\label{eq:apapp1app2}
\begin{array}{lll}
\left(\overline{AP}\right)&\left\{
\begin{array}{l}
\ds  -\partial_x\left(\overline{A_x}\partial_x\overline{u}\right)=\overline{f}+\frac{g_+}{L_z}-\frac{g_-}{L_z}\\
\ds  ~~~~~+\frac{1}{L_z}\partial_x\left(\int_{\Omega_z^1}A_x^\prime \partial_xu_1^\prime \mathrm{d}z+\int_{\Omega_z^2}A_x^\prime \partial_xu_2^\prime\mathrm{d}z\right),\\
\ds \overline{u}\left(x_\pm\right)=0,
\end{array}
\right.
&
\begin{array}{l}
\\
\hspace{-0.5cm}\text{~for~}x\in \Omega_x,\\
~
\end{array}
\\[12mm]
\left(AP_1^\prime \right)&\left\{
\begin{array}{l}
\ds -\partial_x\left(A_x\partial_xu_1^\prime \right)-\partial_z\left(\frac{A_z}{\varepsilon\left(z\right)}\partial_zu_1^\prime \right) =f+\partial_x\left(A_x\partial_x\overline{u}\right),\\
\ds  \frac{A_z\left(x,z_+\right)}{\varepsilon\left(z_+\right)}\partial_zu_1^\prime \left(x,z_+\right)=g_{+}\left(x\right),\\
 u_1^\prime \left(x_\pm,z\right)=0,\\
\ds \frac{A_z\left(x,z_\iota\right)}{\varepsilon\left(z_\iota\right)}\partial_zu_1^\prime\left(x,z_\iota\right)=\frac{A_z\left(x,z_\iota\right)}{\varepsilon\left(z_\iota\right)}\partial_zu_2^\prime\left(x,z_\iota\right),\\
\ds \int_{\Omega_z^1}u_1^\prime\left(x,z\right)\mathrm{d}z+\int_{\Omega_z^2}u_2^\prime\left(x,z\right)\mathrm{d}z=0,
\end{array}
\right.
&
\begin{array}{l}
 \hspace{-0.5cm}\text{~for~}\left(x,z\right)\in \Omega_1,\\
 \hspace{-0.5cm}\text{~for~}x\in \Omega_x,\\
 \hspace{-0.5cm}\text{~for~} z\in \Omega_z^1,\\
 \hspace{-0.5cm}\text{~for~}x\in \Omega_x \text{~(DN-cond)}\\
 \hspace{-0.5cm}\text{~for~}x\in \Omega_x\text{~(constraint)},
\end{array}
\\[22mm]
\left(AP_2^\prime \right)&\left\{
\begin{array}{l}
 \ds -\partial_x\left(A_x\partial_xu_2^\prime \right)-\partial_z\left(\frac{A_z}{\varepsilon\left(z\right)}\partial_zu_2^\prime \right) =f+\partial_x\left(A_x\partial_x\overline{u}\right),\\
 \ds \frac{A_z\left(x,z_-\right)}{\varepsilon\left(z_-\right)}\partial_zu_2^\prime \left(x,z_-\right)=g_{-}\left(x\right),\\
 \ds u_2^\prime \left(x_\pm,z\right)=0,\\
 \ds u_2^\prime\left(x,z_\iota\right)=u_1^\prime\left(x,z_\iota\right),
\end{array}
\right.
&
\begin{array}{l}
\hspace{-0.5cm}\text{~for~}\left(x,z\right)\in \Omega_2,\\
\hspace{-0.5cm}\text{~for~}x\in \Omega_x,\\
\hspace{-0.5cm}\text{~for~} z\in \Omega_z^2,\\
\hspace{-0.5cm}\text{~for~}x\in \Omega_x \text{~(DN-cond)}.
\end{array}
\end{array}
\end{equation}

For the mathematical and numerical study, we shall need the weak form of this system. In this aim, and in order to get rid of the non-homogeneous Dirichlet interface condition, it  will be more convenient to introduce the new variable (lifting procedure)
\begin{equation}
\xi^\prime_2\left(x,z\right):=u_2^\prime\left(x,z\right)-u_1^\prime\left(x,z_\iota\right) \text{~~~for~} \left(x,z\right)\in \Omega_x\times \Omega_z^2,
\label{eq:defofxiprime}
\end{equation}
where we remark that $u_1^\prime\left(x,z_\iota\right) \in H^1(\Omega_x)$, as by Hyp. B, the solution $u$ of the P-model belongs to $H^2(\Omega)$. This leads to the completely equivalent  formulation, called in the following also (AP/AP)-reformulation, for the unknowns $(\bar{u},u_1',\xi_2')$:

\begin{equation}
\label{AP/AP_bis}
\begin{array}{lll}
\left(\overline{AP}\right)&\left\{
\begin{array}{l}
\ds -\partial_x\left(\overline{A_x}\partial_x\overline{u}\right)=\overline{f}+\frac{g_+}{L_z}-\frac{g_-}{L_z}+\frac{1}{L_z}\partial_x\left(\int_{\Omega_z^1}A_x^\prime \partial_xu^\prime_1 \mathrm{d}z\right.\\
\ds ~~~~~\left.+\int_{\Omega_z^2}A_x^\prime \partial_xu_1^\prime\left(x,z_\iota\right) \mathrm{d}z+\int_{\Omega_z^2}A_x^\prime \partial_x\xi_2^\prime\mathrm{d}z\right),\\
\ds \overline{u}\left(x_\pm\right)=0,
\end{array}
\right.
&
\begin{array}{l}
\\
\hspace{-0.8cm}\text{~for~}x\in \Omega_x,\\
~
\end{array}
\\[12mm]
\left(AP_1^\prime \right)&\left\{
\begin{array}{l}
\ds -\partial_x\left(A_x\partial_xu^\prime_1 \right)-\partial_z\left(\frac{A_z}{\varepsilon\left(z\right)}\partial_zu^\prime_1\right) =f+\partial_x\left(A_x\partial_x\overline{u}\right),\\
\ds \frac{A_z\left(x,z_+\right)}{\varepsilon\left(z_+\right)}\partial_zu^\prime_1 \left(x,z_+\right)=g_{+}\left(x\right),\\
\ds u^\prime_1 \left(x_\pm,z\right)=0,\\
\frac{A_z\left(x,z_\iota\right)}{\varepsilon\left(z_\iota\right)}\partial_zu^\prime_1\left(x,z_\iota\right)=\frac{A_z\left(x,z_\iota\right)}{\varepsilon\left(z_\iota\right)}\partial_z\xi^\prime_2\left(x,z_\iota\right),\\
\ds \int_{\Omega_z^1}u^\prime_1\left(x,z\right)\mathrm{d}z+L_z^2u_1^\prime\left(x,z_\iota\right)+\int_{\Omega_z^2}\xi_2^\prime\left(x,z\right)\mathrm{d}z=0,
\end{array}
\right.
&
\begin{array}{l}
\hspace{-0.8cm}\text{~for~}\left(x,z\right)\in \Omega_x\times \Omega_z^1,\\
\hspace{-0.8cm}\text{~for~}x\in \Omega_x,\\
\hspace{-0.8cm}\text{~for~} z\in \Omega_z^1,\\
\hspace{-0.8cm}\text{~for~}x\in \Omega_x\\
\hspace{-0.8cm}\text{~for~}x\in \Omega_x\text{~(constraint)},
\end{array}
\\[21mm]
\left(AP_2^\prime \right)&\left\{
\begin{array}{l}
\ds -\partial_x\left(A_x\partial_x\xi_2^\prime \right)-\partial_z\left(\frac{A_z}{\varepsilon\left(z\right)}\partial_z\xi_2^\prime \right) =f\\
~~~~~+\partial_x\left(A_x\partial_xu_1^\prime\left(x,z_\iota\right) \right)+\partial_x\left(A_x\partial_x\overline{u}\right),\\
\ds \frac{A_z\left(x,z_-\right)}{\varepsilon\left(z_-\right)}\partial_z\xi_2^\prime \left(x,z_-\right)=g_{-}\left(x\right),\\
\ds \xi_2^\prime \left(x_\pm,z\right)=0,\\
\ds \xi_2^\prime\left(x,z_\iota\right)=0,
\end{array}
\right.
&
\begin{array}{l}
~\\
\hspace{-0.8cm}\text{~for~}\left(x,z\right)\in \Omega_x\times \Omega_z^2,\\
\hspace{-0.8cm}\text{~for~}x\in \Omega_x,\\
\hspace{-0.8cm}\text{~for~} z\in \Omega_z^2,\\
\hspace{-0.8cm}\text{~for~}x\in \Omega_x.
\end{array}
\end{array}
\end{equation}

We shall prove in Section \ref{SSEC32} that in the limit of vanishing anisotropy strength in $\Omega_2$, {\it i.e.} $\delta=\eps(z_\iota) \rightarrow 0$, this system yields the following hybrid-system, called (AP/L)-model:
\begin{equation}\label{strongAPL_bis}
\begin{array}{lll}
\left(\overline{AP}\right)&\left\{
\begin{array}{l}
\ds -\partial_x\left(\overline{A_x}\partial_x\overline{u}\right)=\overline{f}+\frac{g_+}{L_z}-\frac{g_-}{L_z}+\frac{1}{L_z}\partial_x\left(\int_{\Omega_z^1}A_x^\prime \partial_xu_1^\prime \mathrm{d}z\right.\\
\ds ~~~~~\left.+\int_{\Omega_z^2}A_x^\prime \partial_x\, u_1^\prime\left(x,z_\iota\right)\mathrm{d}z\right),\\
\ds \overline{u}\left(x_\pm\right)=0,
\end{array}
\right.
&
\begin{array}{l}
\\
\hspace{-0.8cm}\text{~for~}x\in \Omega_x,\\
~
\end{array}
\\[12mm]
\left(AP_1^\prime \right)&\left\{
\begin{array}{l}
\ds -\partial_x\left(A_x\partial_xu_1^\prime \right)-\partial_z\left(\frac{A_z}{\varepsilon\left(z\right)}\partial_zu_1^\prime \right) =f+\partial_x\left(A_x\partial_x\overline{u}\right),\\
\ds \frac{A_z\left(x,z_+\right)}{\varepsilon\left(z_+\right)}\partial_zu_1^\prime \left(x,z_+\right)=g_{+}\left(x\right),\\
\ds u_1^\prime \left(x_\pm,z\right)=0,\\
\ds \frac{A_z\left(x,z_\iota\right)}{\varepsilon\left(z_\iota\right)}\partial_zu_1^\prime\left(x,z_\iota\right)=0,\\
\ds \int_{\Omega_z^1}u_1^\prime\left(x,z\right)\mathrm{d}z+L_z^2\, u_1^\prime\left(x,z_\iota\right)=0,
\end{array}
\right.
&
\begin{array}{l}
\hspace{-0.8cm}\text{~for~}\left(x,z\right)\in \Omega_x\times\Omega_z^1,\\
\hspace{-0.8cm}\text{~for~}x\in \Omega_x,\\
\hspace{-0.8cm}\text{~for~} z\in \Omega_z^1,\\
\hspace{-0.8cm}\text{~for~}x\in \Omega_x\\
\hspace{-0.8cm}\text{~for~}x\in \Omega_x\text{~(constraint)},
\end{array}
\\[21mm]
\left(L\right)&\left\{
\begin{array}{l}
\ds \xi_2^\prime\left(x\right)=0,
\end{array}
\right.
&
\begin{array}{l}
\hspace{-0.8cm}\text{~for~}x\in \Omega_x.
\end{array}
\end{array}
\end{equation}
Formally, one can immediately observe that the solution $\xi_2^\prime$ of $(AP_2')$ converges (in some sense to be defined later on) in the limit $\eps(z_\iota) \rightarrow 0$ towards zero. Indeed, letting formally $\delta=\eps(z_\iota) \rightarrow 0$ implies that the limit function $\xi_2^\prime$ is independent of $z$. The $x$-dependence is then given by the Dirichlet interface condition in $z=z_\iota$, hence $\xi_2^\prime \equiv 0$ in the limit $\eps(z_\iota) \rightarrow 0$.\\
The rigorous mathematical study of the approximation error introduced by using (\ref{strongAPL_bis}) instead of (\ref{AP/AP_bis}) if $\eps(z_\iota) \ll 1$, will be the subject of Section~\ref{SSEC32}, and shall permit to find a criterion determining where to put the interface. Aim of this paper is to show that using the hybrid AP/L-model instead of the fully AP-model is computationally more efficient when the parameter $\eps(z)$ is small in a large part of the domain.

%%%%%%%%%%%%%%%%%%%%%%%%%%%%%%%%%%%%%%%
\subsubsection{Variational formulations}
%%%%%%%%%%%%%%%%%%%%%%%%%%%%%%%%%%%%%%
Our numerical simulations will be based on finite element discretizations of the former systems. For this, we shall introduce now the variational formulations of the (AP/AP)-model (\ref{AP/AP_bis}) as well as of the (AP/L)-model (\ref{strongAPL_bis}), and define the Hilbert-spaces
\begin{eqnarray*}
\mathcal{V}_1&:=&\left\{\psi\left(\cdot,\cdot\right)\in H^1\left(\Omega_1\right)~/~\psi(x_\pm,\cdot)=0\right\}\,,\\% \quad \left(\phi,\psi\right)_{\mathcal{V}_1}:=\left(\partial_x\phi,\partial_x\psi\right)_{L^2_1}+\left(\partial_z\phi,\partial_z\psi\right)_{L^2_1},\\
\mathcal{V}_2&:=&\left\{\psi\left(\cdot,\cdot\right)\in H^1\left(\Omega_2\right)~/~\psi(x_\pm,\cdot)=0\,, \,\,\, \psi(\cdot,z_{\iota})=0 \right\}\,,\\%\quad \left(\phi,\psi\right)_{\mathcal{V}_2}:=\left(\partial_x\phi,\partial_x\psi\right)_{L^2_2}+\left(\partial_z\phi,\partial_z\psi\right)_{L^2_2},\\
\mathcal{W}&:=&\left\{\psi\left(\cdot\right)\in H^1\left(\Omega_x\right)~/~\psi(x_\pm)=0 \right\}\,, %\quad \left(\phi,\psi\right)_\mathcal{W}:=\left(\partial_x\phi,\partial_x\psi \right)_{L^2_x}
\end{eqnarray*}
associated to the scalar products
\begin{eqnarray*}
\left(\phi,\psi\right)_{\mathcal{V}_1}&:=&\left(\partial_x\phi,\partial_x\psi\right)_{L^2_1}+\left(\partial_z\phi,\partial_z\psi\right)_{L^2_1},\\
\left(\phi,\psi\right)_{\mathcal{V}_2}&:=&\left(\partial_x\phi,\partial_x\psi\right)_{L^2_2}+\left(\partial_z\phi,\partial_z\psi\right)_{L^2_2},\\
\left(\phi,\psi\right)_\mathcal{W}&:=&\left(\partial_x\phi,\partial_x\psi\right)_{L^2_x},
\end{eqnarray*}
where we denoted by $(\cdot,\cdot)_{L^2_1}$ resp. $(\cdot,\cdot)_{L^2_2}$ the $L^2$ scalar-products in $\Omega_1$ resp. $\Omega_2$. To simplify the writing, we define further for $v^\prime\in\mathcal{V}_\star$, $\psi^\prime\in\mathcal{V}_\star$, $\overline{v}\in\mathcal{W}$, $\overline{\psi}\in\mathcal{W}$, $\overline{P}\in L^2\left(\Omega_x\right)$ and $\overline{Q}\in L^2\left(\Omega_x\right)$, the following bilinear forms
\begin{equation}\label{eq:bilinear:forms}
  \begin{array}{llll}
a_{z\star}\left(v^\prime,\psi^\prime\right)&:=\int_{\Omega_z^\star}\int_{\Omega_x}\frac{A_z}{\varepsilon\left(z\right)}\partial_zv^\prime\partial_z\psi^\prime \mathrm{d}x\mathrm{d}z,&&\\[2mm]
a_{xf\star}\left(v^\prime,\psi^\prime\right)&:=\int_{\Omega_z^\star}\int_{\Omega_x}A_x\partial_xv^\prime\partial_x\psi^\prime \mathrm{d}x\mathrm{d}z,& a_{xa}\left(\overline{v},\overline{\psi}\right)&:=\int_{\Omega_x}\overline{A_x}\partial_x\overline{v}\partial_x\overline{\psi}\mathrm{d}x,\\[2mm]
b_{l1}\left(\overline{P},\psi^\prime\right)&:=\int_{\Omega_x}\overline{P}\int_{\Omega_z^1}\frac{1}{\varepsilon\left(z\right)}\psi^\prime \mathrm{d}z\mathrm{d}x,& b_{c\star}\left(v^\prime,\overline{Q}\right)&:=\frac{1}{L_z}\int_{\Omega_x}\overline{Q}\int_{\Omega_z^\star}v^\prime \mathrm{d}z\mathrm{d}x,\\[2mm]
%b_{C2}\left(u_2^\prime,\overline{Q}\right)&:=&\frac{L_z^2}{L_z}\int_{\Omega_x}\overline{Q}u_2^\prime ~\mathrm{d}x,\\
c_{f\star}\left(\overline{v},\psi^\prime\right)&:=\int_{\Omega_z^\star}\int_{\Omega_x}A_x\partial_x\overline{v}\partial_x\psi^\prime \mathrm{d}x\mathrm{d}z,&c_{a\star}\left(v^\prime,\overline{\psi}\right)&:=\int_{\Omega_z^\star}\int_{\Omega_x}A_x^\prime\partial_xv^\prime\partial_x\overline{\psi}\mathrm{d}x\mathrm{d}z,\\[2mm]
d_\iota\left(v^\prime,\psi_1^\prime\right)&:=\!\!\! \int_{\Omega_x}\frac{A_z\left(x,z_\iota\right)}{\varepsilon\left(z_\iota\right)}\partial_zv^\prime\left(x,z_\iota\right)\psi_1^\prime\left(x,z_\iota\right)\mathrm{d}x,\quad&  v^\prime \,\,\, s.t.\,\,\,& \frac{A_z\left(\cdot,z_\iota\right)}{\varepsilon\left(z_\iota\right)}\partial_zv^\prime\left(\cdot,z_\iota\right) \in L^2(\Omega_x)\,,  
  \end{array}
\end{equation}
where $\star$ takes the value $1$ or $2$. 
The variational formulation of the (AP/AP)-model (\ref{AP/AP_bis}) writes: \,\,\, \textit{Find $(\overline{u},u^\prime_1,\xi_2^\prime,\overline{P})  \in\mathcal{W} \times \mathcal{V}_1 \times \mathcal{V}_2 \times L^2\left(\Omega_x\right)$ such that}

\begin{equation}
\label{AP/AP-FV}
\left(AP/AP\right)\left\{
\begin{array}{l}
a_{xa}\left(\overline{u},\overline{\psi}\right)=\left(\overline{f},\overline{\psi}\right)_{L^2_x}+\frac{1}{L_z}\left(g_+-g_-,\overline{\psi}\right)_{L^2_x}\\[3mm]
~~~~~~~-\frac{1}{L_z}\left(c_{a1}\left(u_1^\prime,\overline{\psi}\right)+c_{a2}\left(u_1^\prime\left(\cdot,z_\iota\right),\overline{\psi}\right)+c_{a2}\left(\xi_2^\prime,\overline{\psi}\right)\right), \qquad \forall \overline{\psi}\in\mathcal{W},\\[3mm]
a_{xf1}\left(u_1^\prime,\psi_1^\prime\right)+a_{z1}\left(u_1^\prime,\psi_1^\prime\right)+b_{l1}\left(\overline{P},\psi_1^\prime\right)=\left(f,\psi_1^\prime\right)_{L^2_1}\\[3mm]
~~~~~~~+\left(g_+,\psi_1^\prime\left(\cdot,z_+\right)\right)_{L^2_x}-d_\iota\left(\xi_2^\prime,\psi_1^\prime\right)-c_{f1}\left(\overline{u},\psi_1^\prime\right), \qquad \forall \psi_1^\prime\in\mathcal{V}_1,\\[3mm]
a_{xf2}\left(\xi_2^\prime,\psi_2^\prime\right)+a_{z2}\left(\xi_2^\prime,\psi_2^\prime\right)=\left(f,\psi_2^\prime\right)_{L^2_2}-\left(g_-,\psi_2^\prime\left(\cdot,z_-\right)\right)_{L^2_x}-c_{f2}\left(\overline{u},\psi_2^\prime\right)\\[3mm]
~~~~~~~-a_{xf2}\left(u_1^\prime\left(\cdot,z_\iota\right),\psi_2^\prime\right), \qquad \forall \psi_2^\prime\in\mathcal{V}_2,\\[3mm]
b_{c1}\left(u_1^\prime,\overline{Q}\right)=-b_{c2}\left(u_1^\prime\left(\cdot,z_\iota\right),\overline{Q}\right)-b_{c2}\left(\xi_2^\prime,\overline{Q}\right), \qquad \forall \overline{Q}\in L^2.\end{array}
\right.
\end{equation}

Equally, we obtain the following variational formulation for the (AP/L)-model (\ref{strongAPL_bis}): 

\textit{Find $(\overline{u},u^\prime_1,\overline{P}) \in\mathcal{W} \times \mathcal{V}_1\times L^2\left(\Omega_x\right)$ such that}
\begin{equation}\left(AP/L\right)
\left\{
\begin{array}{ll}
a_{xa}\left(\overline{u},\overline{\psi}\right)=\left(\overline{f},\overline{\psi}\right)_{L^2_x}+\frac{1}{L_z}\left(g_+-g_-,\overline{\psi}\right)_{L^2_x}\\[3mm]
~~~~~-\frac{1}{L_z}\left(c_{a1}\left(u^\prime_1,\overline{\psi}\right)+c_{a2}\left(u^\prime_1\left(\cdot,z_\iota\right),\overline{\psi}\right)\right), &\forall \overline{\psi}\in\mathcal{W},\\[3mm]
a_{xf1}\left(u^\prime_1,\psi_1^\prime\right)+a_{z1}\left(u^\prime_1,\psi_1^\prime\right)+b_{l1}\left(\overline{P},\psi_1^\prime\right)\\[3mm]
~~~~~=\left(f,\psi_1^\prime\right)_{L^2_1}+\left(g_+,\psi_1^\prime\left(\cdot,z_+\right)\right)_{L^2_x}-c_{f1}\left(\overline{u},\psi_1^\prime\right), &\forall \psi_1^\prime\in\mathcal{V}_1,\\[3mm]
b_{c1}\left(u^\prime_1,\overline{Q}\right)=-b_{c2}\left(u^\prime_1\left(\cdot,z_\iota\right),\overline{Q}\right), &\forall \overline{Q}\in L^2\left(\Omega_x\right).
\end{array}
\right.
\label{weakform}
\end{equation}

It will be this (AP/L)-model which will be used in our numerical simulations (Section~\ref{SEC4}) in order to compute the solution of the highly anisotropic elliptic equation (\ref{eq:SP}).

%%%%%%%%%%%%%%%%%%%%%%%%%%%%%%%%%%%%%%%%%%%%%%%%%%%
\subsection{Mathematical study}\label{SSEC32}
%%%%%%%%%%%%%%%%%%%%%%%%%%%%%%%%%%%%%%%%%%%%%%%%%%%

Aim of this section is to study the approximation error introduced when using the Limit-model in the sub-domain $\Omega_x \times \Omega_z^2$ (where the anisotropy strength $\eps$ is small) instead of the original problem. In other words we are interested in the error between the solution of the (AP/AP)-model and the (AP/L)-model. This error will depend on the position of the coupling interface, more precisely on $\eps\left(z_\iota\right)$.

Let us denote in the following by $\left(\bar{u},u_1^\prime,\xi_2^\prime\right)\in {\mathcal W} \times{\mathcal V}_1 \times{\mathcal V}_2$ the solution of (\ref{AP/AP_bis}) or \eqref{AP/AP-FV} and by $\left({\bar v},v_1^\prime,0\right)\in {\mathcal W} \times{\mathcal V}_1 \times{\mathcal V}_2$ the solution of (\ref{strongAPL_bis}) or \eqref{weakform}. Remark that the existence and uniqueness of a solution of (\ref{AP/AP_bis}) can be straightforwardly shown by equivalence with the AP-model (\ref{AP_bar})-(\ref{AP_fluc}). Moreover the existence and uniqueness of a solution of (\ref{strongAPL_bis}) is shown by standard arguments, proving the boundedness of $\left(\bar{u},u_1^\prime,\xi_2^\prime\right)$ in ${\mathcal W} \times{\mathcal V}_1 \times{\mathcal V}_2$, uniformly in $\eps(z_\iota)$, inducing thus weak convergences and passing finally to the limit in the variational formulation. We leave the details to the reader and prefer to concentrate on the error estimate. We will begin with estimating firstly  $\xi_2^\prime$, showing its convergence towards zero  for $\eps\left(z_\iota\right) \rightarrow 0$ (Theorem \ref{th_xi2}). Using this convergence result, we shall estimate in a second step  the errors $||\overline{u}-\overline{v}||_{\mathcal W}$ as well as $|| u_1^\prime-v_1^\prime ||_{{\mathcal V}_1} $ (Theorem \ref{mean}). For the convergence of  $\xi_2^\prime$ we will need the following regularity result.

\begin{lemme} \label{lemme_reg}
Under Hypothesis B, the solution $u$ of the P-model (\ref{P}) belongs to $H^2\left(\Omega\right)$ and one has the following estimates
\begin{equation}
\begin{split}
||u||_{H^2\left(\Omega\right)}+||u^\prime||_{H^2\left(\Omega\right)}+||\overline{u}||_{H^2\left(\Omega_x\right)}+||u^\prime\left(\cdot,z_\iota\right)||_{H^1\left(\Omega_x\right)} \le C_{f,g_\pm} (1+ \eps_M)\,,
%\le c\left(||f||_{L^2\left(\Omega\right)}+||g_-||_{L^2\left(\Omega_x\right)}+||g_+||_{L^2\left(\Omega_x\right)}\right),
\end{split}
\end{equation}
with some constant $C_{f,g_\pm}>0$ independent of $\eps$, however dependent on the source terms $f$ respectively $g_\pm$. Recall that $\overline{u}$ is the average of $u$ over $z$ and $u^\prime:=u-\overline{u}$.
\end{lemme}

\debproof
The proof of this lemma uses standard ``energy''-estimates and density arguments as well as the Poincar\'e and trace-theorems. See \cite{brezis,evans,DLNN} for more details.\\

Standard elliptic results permit to show that, under the additional hypothesis of this  theorem, the solution of (\ref{P}) is more regular, {\it i.e.} ${u} \in H^2(\Omega)$.

 Remark then that multiplying (\ref{P}) by ${u}$ and integrating over
  $\Omega$ yields by integration by parts the
  $H^1$-estimate
  \begin{equation} \label{H1}
    ||\partial_x {u}||^2_{L^{2}}+ || {1 \over \sqrt{\eps(z)}}\, \partial_z {u}||^2_{L^{2}}  \le C ( ||f||^2_{L^{2}}+ ||g_+||^2_{L^2_x}+||g_-||^2_{L^2_x})\,.
  \end{equation}
  Rewriting now the equation as 
  $$
  - A_x \partial_{xx} {u} - \partial_z\left(  {A_z \over \eps(z)} \partial_{z} {u} \right) =f+  (\partial_x A_x) \partial_x {u} \,,
  $$
  multiplying it by $-\partial_{xx} {u}$ and
  integrating over $\Omega$ yields (by integration by parts)
  \be \label{ESTi}
  \begin{array}{l}
 \ds   || \sqrt{A_x}\, \partial_{xx} {u}||^2 +|| \sqrt{A_z \over \eps(z)}\, \partial_{xz} {u}||^2 + \int_{\Omega} {\partial_x A_z \over \eps(z)} \, \partial_z u \, \partial_{xz}u\, dxdz \\[3mm]
\ds     \hspace{1cm} \le ||f|| \,||\partial_{xx} {u}||  +  C ||\partial_x u||\, ||\partial_{xx}u|| + \int_{\Omega_x}\!\! \partial_x g_+ \partial_x u(x,z_+)\, dx -\int_{\Omega_x}\!\! \partial_x g_- \partial_x u(x,z_-)\, dx \,.
  \end{array} 
  \ee
  Using now the $H^1$-estimate (\ref{H1}) and the fact that $||\partial_x u (x,z_{\pm})||_{L^2_x} \le C||\partial_x u||_{L^2} + C||\partial_{xz} u||_{L^2}$, yields the first $H^2$-estimates
\be \label{H2_est1}
||\partial_{xx} {u}||_{L^2}^2 + || { 1 \over \sqrt{\eps(z)} }\,\partial_{xz} {u}||_{L^2}^2 \le C_{f,g_{\pm}} (1 + \eps_{max})\,.
\ee
Coming now back to the equation
  $$
  - \partial_z ({A_z\over \eps(z)}\, \partial_z {u})  =f + \partial_x (A_x \,\partial_x {u})\,,
  $$
  one gets with \eqref{H2_est1}
  $$
  ||\partial_z ({A_z \over \eps(z)}\, \partial_z {u})|| \le C_{f,g_{\pm}}(1 + \eps_{max})\,.
  $$
Hence, remarking that 
$$
{A_z \over \eps(z)}\, \partial_z {u} (x,z)={A_z \over \eps(z_-)}\, \partial_z {u} (x,z_-)+ \int_{z_-}^z \partial_z \left( {A_z \over \eps(\tau )}\, \partial_z {u} (x,\tau)\right)\, d\tau\,,
$$
permits to show an improved $H^1$-estimate
\be \label{EEE}
  ||{ 1 \over \eps(z)} \,\partial_z {u} ||_{L^2} \le C ||g_-||_{L^2_x}+C ||\partial_z ( {A_z \over \eps(z)}\, \partial_z {u})||_{L^2} \le C_{f,g_{\pm}}(1 + \eps_{max})\,.
\ee
Finally, with all these estimates, one has
$$
\partial_{zz}u = \partial_z \left( {\eps(z) \over A_z} \, {A_z \over \eps(z)}\,  \partial_z u\right)= \partial_z \left({\eps(z) \over A_z} \right) {A_z \over \eps(z)}\,  \partial_z {u} + {\eps(z) \over A_z}\, \partial_z \left({A_z \over \eps(z)}\, \partial_z {u}\right)\,,
$$
yielding 
\be \label{H2-est2}
||\partial_{zz}u||_{L^2} \le ||\partial_z \left({\eps(z) \over A_z} \right)||_{L^\infty}\, ||{A_z \over \eps(z)} \partial_z {u} ||_{L^2} + C \eps_{max} \, ||\partial_z \left({A_z \over \eps(z)} \partial_z {u}\right)||_{L^2} \le C_{f,g_{\pm}}\eps_{M}\,,
\ee
which concludes the proof together with \eqref{H2_est1}.

\finproof

~

From now on, we shall not carry with us the bound of $||\eps||_{W^{1,\infty}(\Omega_z)}$, {\it i.e.} $\eps_M$, keeping in mind that the anisotropy has to be bounded from above in $W^{1,\infty}$.

\begin{theorem} Under Hypothesis B, the function $\xi_2' \in {\mathcal V}_2$, part of the solution of the (AP/AP)-model (\ref{AP/AP_bis}), converges towards zero in $H^1\left(\Omega_2\right)$ as $\varepsilon\left(z_\iota\right)\rightarrow 0$. If we assume moreover that $u_1^\prime\left(\cdot,z_\iota\right)\in H^2\left(\Omega_x\right)$, then there exists a constant $c>0$, independent of $\eps(z_\iota)$, such that
\begin{equation}
\label{eq:est_on_xi2}
||\partial_x\xi_2^\prime||_{L^2\left(\Omega_2\right)}   \leq c\sqrt{\varepsilon\left(z_\iota\right)}~~~\text{and}~~~||\partial_z\xi_2^\prime||_{L^2\left(\Omega_2\right)} \leq c\, \varepsilon\left(z_\iota\right).
\end{equation}
%meaning $\xi_2^\prime\xrightarrow[\varepsilon\left(z_\iota\right)\rightarrow 0]{} 0$ in $H^1\left(\Omega\right)$.
\label{th_xi2}
\end{theorem}

\debproof The proof of this theorem is very similar to the proof of the previous lemma. Let us only recall that $\xi_2'$ is the solution of a diffusion problem with slightly different boundary condition on one side, when compared with the P-problem, {\it i.e.}
\be \label{XI}
\left\{
\begin{array}{ll}
\ds -\partial_x\left(A_x\, \partial_x\xi_2^\prime \right)-\partial_z\left(\frac{A_z}{\varepsilon\left(z\right)}\, \partial_z\xi_2^\prime \right) =h& \text{~for~}\left(x,z\right)\in \Omega_x\times \Omega_z^2\\
\ds \frac{A_z\left(x,z_-\right)}{\varepsilon\left(z_-\right)}\, \partial_z\xi_2^\prime \left(x,z_-\right)=g_{-}\left(x\right),&\text{~for~}x\in \Omega_x,\\
\ds \xi_2^\prime\left(x,z_\iota\right)=0,&\text{~for~}x\in \Omega_x\\
\ds \xi_2^\prime \left(x_\pm,z\right)=0,&\text{~for~} z\in \Omega_z^2\,,
\end{array}
\right.
\ee
where we denoted for simplicity reasons $h:=f+\partial_x\left(A_x\partial_xu_1^\prime\left(x,z_\iota\right) \right)+\partial_x\left(A_x\partial_x\overline{u}\right)$. Let us first assume that $u_1^\prime\left(\cdot,z_\iota\right)\in H^2\left(\Omega_x\right)$, such that $h \in L^2(\Omega_2)$. Multiplying now the system \eqref{XI} by $\xi_2^\prime$ and integrating over $\Omega_x\times\Omega_z^2$ yields
$$
|| \sqrt{A_x}\, \partial_x\xi_2^\prime ||^2_{L^2_2} + || \sqrt{\frac{A_z}{\varepsilon\left(z\right)}} \, \partial_z\xi_2^\prime||^2_{L^2_2}= \int_{\Omega^2}h\, \xi_2^\prime\mathrm{d}x\mathrm{d}z-\int_{\Omega_x}g_-\left(x\right)\xi_2^\prime\left(x,z_-\right)\, \mathrm{d}x\,.
$$
Remarking that $\frac{1}{\varepsilon\left(z\right)}\geq\frac{1}{\varepsilon\left(z_\iota\right)}$ in $\Omega_z^2$, and that $\xi_2^\prime\left(x,z_\iota\right)=0$ for all $x\in\Omega_x$, such that by Poincar\'e's inequality
\begin{equation} \label{Poi}
\|\xi_2^\prime\left(\cdot,z_-\right)\|_{L^2\left(\Omega_x\right)}\leq c \|\partial_z\xi_2^\prime\|_{L^2\left(\Omega_2\right)}~~~\text{and}~~~\|\xi_2^\prime\|_{L^2\left(\Omega_2\right)}\leq c \|\partial_z\xi_2^\prime\|_{L^2\left(\Omega_2\right)},
\end{equation}
we obtain by Cauchy-Schwarz
\begin{equation*}
\begin{split}
\|\partial_x\xi_2^\prime\|^2_{L^2\left(\Omega_2\right)}+\frac{1}{\varepsilon\left(z_\iota\right)}\|\partial_z\xi_2^\prime\|^2_{L^2\left(\Omega_2\right)}\leq c \big[ \| h\|_{L^2\left(\Omega_2\right)}  + \| g_-\|_{L^2\left(\Omega_x\right)} \big]\|\partial_z\xi_2^\prime\|_{L^2\left(\Omega_2\right)},
\end{split}
\end{equation*}
leading to 
\begin{equation} \label{estimate}
\varepsilon\left(z_\iota\right)\|\partial_x\xi_2^\prime\|^2_{L^2\left(\Omega_2\right)}+\|\partial_z\xi_2^\prime\|^2_{L^2\left(\Omega_2\right)}\leq c(h,g_-)\varepsilon\left(z_\iota\right)^2,
\end{equation}
which completes the proof in the regular case. Remark that it is at this point that our proof differs from the $H^1$-estimates \eqref{H1}, \eqref{EEE} of the previous proof, and this thanks to the Poincar\'e inequality \eqref{Poi}, valid due to the boundary condition 
$\xi_2^\prime\left(x,z_\iota\right)=0$, $\forall x\in\Omega_x$. It is exactly this condition which implies the convergence of $\xi_2'$ towards zero.\\

Let us suppose now that we have only  $u_1^\prime\left(\cdot,z_\iota\right)\in H^1\left(\Omega_x\right)$, which is immediate under Hypothesis B, such that 
$h$ belongs only to $H^{-1}(\Omega_2)$. In this case, following the same arguments as above, would lead rather to
\begin{equation*}
\begin{split}
\|\partial_x\xi_2^\prime\|^2_{L^2\left(\Omega_2\right)}+\frac{1}{\varepsilon\left(z_\iota\right)}\|\partial_z\xi_2^\prime\|^2_{L^2\left(\Omega_2\right)}\leq c \| h\|_{H^{-1}(\Omega_2)}\|\xi_2^\prime\|_{H^1\left(\Omega_2\right)},  + c \| g_-\|_{L^2\left(\Omega_x\right)} \|\partial_z\xi_2^\prime\|_{L^2\left(\Omega_2\right)},
\end{split}
\end{equation*}
leading only to 
\begin{equation} \label{estimate_bis}
\|\partial_x\xi_2^\prime\|^2_{L^2\left(\Omega_2\right)}+ { 1 \over \varepsilon\left(z_\iota \right)} \|\partial_z\xi_2^\prime\|^2_{L^2\left(\Omega_2\right)}\leq c (\| h\|^2_{H^{-1}\left(\Omega_2\right)} +  \sqrt{\varepsilon\left(z_\iota\right)} \| g_-\|^2_{L^2\left(\Omega_x\right)} ) ,
\end{equation}
which does not permit to have the desired convergence. However, density arguments permit to conclude the proof in the general case. Indeed, let us consider the linear, continuous mapping
 \be \label{MU}
\begin{array}{c}
\displaystyle {\mathcal U}: h \in H^{-1}(\Omega_2)  \longmapsto  \xi_2^\prime \in H^1(\Omega_2)\,\, \,\, \textrm{sol. of } \eqref{XI}\,.
\end{array}
\ee
Let now $h \in  H^{-1}(\Omega_2)$. Then for each $\delta >0$ there exists by density a more regular data $\tilde{h} \in L^2(\Omega_2)$ such that $||h-\tilde{h}||_{ H^{-1}(\Omega_2)} \le \delta$ and $||{\mathcal U}(\tilde{h}) ||_{H^1(\Omega_2)} \le c(\tilde{h}) \, \sqrt{\eps(z_\iota)}$ by the estimate \eqref{estimate}. As the map ${\mathcal U}$ is continuous and linear, one has then
%For more regular data $F$, we can have $u'(\cdot,z_\iota) \in H^2(\Omega_x)$, yielding the estimate previously proven. Let now $F$ having minimum regularity (satisfying Hyp. B). Then, for each $\delta >0$ there exists $\tilde{F}$ being regular enough such that $||{\mathcal U}(\tilde{F}) ||_{H^1(\Omega_2)} \le c \, \sqrt{\eps(z_\iota)}$ and $||F-\tilde{F}|| \le \delta$, where the norm corresponds here to the composition of the norms of the data-set. As the map ${\mathcal U}$ is continuous and linear, one has then
$$
\begin{array}{lll}
\ds ||{\mathcal U}(h)||_{H^1(\Omega_2)} \le ||{\mathcal U}(h-\tilde{h})||_{H^1(\Omega_2)}+||{\mathcal U}(\tilde{h})||_{H^1(\Omega_2)} &\le&\ds  c||h-\tilde{h}||_{H^{-1}(\Omega_2)} + c(\tilde{h}) \sqrt{\eps(z_\iota)}\\[3mm]
& \le&\ds  c \delta + c(\tilde{h}) \sqrt{\eps(z_\iota)}\,.
\end{array}
$$ 
Thus, for all $\tau >0$, one can find $\delta >0$ and $\eps(z_\iota) >0$ such that $||{\mathcal U}(h)||_{H^1(\Omega_2)} \le \tau $ for all $\eps < \eps(z_\iota)$, which proves the convergence of ${\mathcal U}(h)$ towards zero in ${H^1(\Omega_2)}$, as $\eps(z_\iota) \rightarrow 0$.

\finproof

Let us now come to the study of the error on the mean part $\overline{u}-\overline{v}$ and the fluctuation part $u_1^\prime-v_1^\prime$, when approximating the (AP/AP)-model with the (AP/L)-model.

\begin{theorem}\label{TheTheorem}
Let $\left(\bar{u},u_1^\prime,\xi_2^\prime\right)\in {\mathcal W} \times{\mathcal V}_1 \times{\mathcal V}_2$ be the solution of (\ref{AP/AP_bis})  and $\left({\bar v},v_1^\prime,0\right)\in {\mathcal W} \times{\mathcal V}_1 \times{\mathcal V}_2$ the solution of (\ref{strongAPL_bis}) and let us assume Hypothesis B to be satisfied. Then one has
$$
\bar{u} \rightarrow_{\eps(z_\iota) \rightarrow 0} \bar v \,\,\, \textrm{in}\,\,\, {\mathcal W} \quad \textrm{and} \quad u_1^\prime \rightarrow_{\eps(z_\iota) \rightarrow 0} v_1^\prime \,\,\, \textrm{in}\,\,\, {\mathcal V}_1\,.
$$
If we suppose moreover that $u_1^\prime\left(\cdot,z_\iota\right)\in H^2\left(\Omega_x\right)$, then there exists a positive constant $c$, independent of $\eps(z_\iota)$, such that
\begin{equation} \label{ESS}
\|\partial_x\, (\overline{u}-\overline{v})\|_{L^2\left(\Omega_x\right)}+ \|\partial_x\, (u_1^\prime-v_1^\prime)\|_{L^2\left(\Omega_1\right)}+ \|\partial_z\, (u_1^\prime-v_1^\prime)\|_{L^2\left(\Omega_1\right)} \leq c\sqrt{\varepsilon\left(z_\iota\right)}\,.
\end{equation}
 \label{mean}
\end{theorem}

\debproof 
To simplify the notations, we shall denote in this proof the difference between the two solutions of (\ref{AP/AP_bis}) resp. (\ref{strongAPL_bis}) by $\left(\overline{w}, w_1^\prime, \xi_2^\prime\right)$, where $\overline{w}:=\overline{u}-\overline{v} \in {\mathcal W}$ and $w_1^\prime:=u_1^\prime-v_1^\prime \in {\mathcal V}_1$. Aim is to show that both $\bar{w}$ and $w_1^\prime$ converge towards zero in the respective spaces. This shall be done in several steps.

First, let us suppose that $u_1^\prime\left(\cdot,z_\iota\right)\in H^2\left(\Omega_x\right)$ and start by writing the system of equations satisfied by $\overline{w}$ and $w_1^\prime$:
\begin{equation}
\label{AP/APL_1}
\left\{
\begin{array}{l}
\ds -\partial_x\left(\overline{A_x}\partial_x\overline{w}\right)=\frac{1}{L_z}\partial_x\left(\int_{\Omega_z^1}A_x^\prime \partial_xw^\prime_1 \mathrm{d}z\right.\\[2mm]
\ds ~~~~~\left.+\int_{\Omega_z^2}A_x^\prime \partial_xw_1^\prime\left(x,z_\iota\right) \mathrm{d}z+\int_{\Omega_z^2}A_x^\prime \partial_x\xi_2^\prime\mathrm{d}z\right),\quad \text{~for~}x\in \Omega_x,\\[2mm]
\ds \overline{w}\left(x_\pm\right)=0,
\end{array}
\right.
\end{equation}
\begin{equation}
\label{AP/APL_2}
\left\{
\begin{array}{ll}
\ds -\partial_x\left(A_x\partial_xw^\prime_1 \right)-\partial_z\left(\frac{A_z}{\varepsilon\left(z\right)}\partial_zw^\prime_1\right) =\partial_x\left(A_x\partial_x\overline{w}\right),& \text{~for~}\left(x,z\right)\in \Omega_x\times \Omega_z^1,\\[2mm]
\ds \frac{A_z\left(x,z_+\right)}{\varepsilon\left(z_+\right)}\partial_zw^\prime_1 \left(x,z_+\right)=0,&\text{~for~}x\in \Omega_x,\\[2mm]
w^\prime_1 \left(x_\pm,z\right)=0,&\text{~for~} z\in \Omega_z^1,\\[2mm]
\ds \frac{A_z\left(x,z_\iota\right)}{\varepsilon\left(z_\iota\right)}\partial_zw^\prime_1\left(x,z_\iota\right)=\frac{A_z\left(x,z_\iota\right)}{\varepsilon\left(z_\iota\right)}\partial_z\xi^\prime_2\left(x,z_\iota\right),&\text{~for~}x\in \Omega_x \\[3mm]
\ds \int_{\Omega_z^1}w^\prime_1\left(x,z\right)\mathrm{d}z+L_z^2w_1^\prime\left(x,z_\iota\right)+\int_{\Omega_z^2}\xi_2^\prime\left(x,z\right)\mathrm{d}z=0,&\text{~for~}x\in \Omega_x\text{~(constraint)}.
\end{array}
\right.
\end{equation}
Introducing on $\Omega_x \times \Omega_z^1$ the variable $w(x,z):=\overline{w}(x)+w_1^\prime(x,z)$, this one solves:
\begin{equation} \label{WW}
\begin{array}{ll}
\left\{
\begin{array}{l}
-\partial_x\left(A_x\partial_xw \right)-\partial_z\left(\frac{A_z}{\varepsilon\left(z\right)}\partial_zw\right) =0,\\
\frac{A_z\left(x,z_+\right)}{\varepsilon\left(z_+\right)}\partial_zw\left(x,z_+\right)=0,\\
w \left(x_\pm,z\right)=0,\\
\frac{A_z\left(x,z_\iota\right)}{\varepsilon\left(z_\iota\right)}\partial_zw\left(x,z_\iota\right)=\frac{A_z\left(x,z_\iota\right)}{\varepsilon\left(z_\iota\right)}\partial_z\xi^\prime_2\left(x,z_\iota\right),
\end{array}
\right.
&
\begin{array}{l}
\text{~for~}\left(x,z\right)\in \Omega_x\times \Omega_z^1,\\
\text{~for~}x\in \Omega_x,\\
\text{~for~} z\in \Omega_z^1,\\
\text{~for~}x\in \Omega_x.
\end{array}
\end{array}
\end{equation}
{\bf 1st Step}: $H^1$-estimate of $w$:\\
Multiplying \eqref{WW} by $w$ and integrating on $\Omega_x\times\Omega_z^1$ by parts, yields
\begin{equation}
\begin{split}
|| \sqrt{A_x}\, \partial_x w ||^2_{L^2_1}+ || \sqrt{A_z \over \varepsilon\left(z\right) }\, \partial_z w ||^2_{L^2_1}
+\int_{\Omega_x}\frac{A_z\left(x,z_\iota\right)}{\varepsilon\left(z_\iota\right)}\partial_z\xi_2^\prime\left(x,z_\iota\right)w\left(x,z_\iota\right)\mathrm{d}x=0.
\end{split}\label{eq:esti1_on_w}
\end{equation}
The difficulty in estimating $w$ comes now from the term $\int_{\Omega_x}\frac{A_z\left(x,z_\iota\right)}{\varepsilon\left(z_\iota\right)}\partial_z\xi_2^\prime\left(x,z_\iota\right)w\left(x,z_\iota\right)\mathrm{d}x$. To be able to reformulate this term,  we shall come back to \eqref{AP/APL_1}-\eqref{AP/APL_2}.

~

The constraint in (\ref{AP/APL_2}) can be reformulated as
\begin{equation*}
\begin{split}
0=\partial_x\left(\overline{A_x} \partial_x\left[\int_{\Omega_z^1}w^\prime_1 \mathrm{d}z+\int_{\Omega_z^2}w_1^\prime\left(x,z_\iota\right) \mathrm{d}z+\int_{\Omega_z^2}\xi_2^\prime\mathrm{d}z\right]\right)\\[3mm]
=\partial_x\left(\int_{\Omega_z^1}\overline{A_x} \partial_xw^\prime_1 \mathrm{d}z+\int_{\Omega_z^2}\overline{A_x} \partial_xw_1^\prime\left(x,z_\iota\right) \mathrm{d}z+\int_{\Omega_z^2}\overline{A_x} \partial_x\xi_2^\prime\mathrm{d}z\right)\,,
\end{split}
\end{equation*}
which allows to rewrite \eqref{AP/APL_1} as
\begin{equation}
\label{eq:onwbar}
-\partial_x\left(\overline{A_x}\partial_x\overline{w}\right)=\frac{1}{L_z}\partial_x\left(\int_{\Omega_z^1}A_x \partial_xw^\prime_1 \mathrm{d}z+\int_{\Omega_z^2}A_x \partial_xw_1^\prime\left(x,z_\iota\right) \mathrm{d}z+\int_{\Omega_z^2}A_x \partial_x\xi_2^\prime\mathrm{d}z\right).
\end{equation}
Integrating now \eqref{AP/APL_2} on $\Omega_z^1$ gives
\begin{equation}
-\int_{\Omega_z^1}\partial_x\left(A_x\partial_xw^\prime_1 \right)\mathrm{d}z+\frac{A_z\left(x,z_\iota\right)}{\varepsilon\left(z_\iota\right)}\partial_z\xi^\prime_2\left(x,z_\iota\right) =\int_{\Omega_z^1}\partial_x\left(A_x\partial_x\overline{w}\right)\mathrm{d}z,
\end{equation}
and (\ref{eq:onwbar}) becomes then
\begin{equation}
\begin{split}
-\frac{1}{L_z}\int_{\Omega_z}\partial_x\left(A_x\partial_x\overline{w}\right)\mathrm{d}z=\frac{1}{L_z}\left(-\int_{\Omega_z^1}\partial_x\left(A_x\partial_x\overline{w}\right)\mathrm{d}z+\frac{A_z\left(x,z_\iota\right)}{\varepsilon\left(z_\iota\right)}\partial_z\xi^\prime_2\left(x,z_\iota\right)\right.\\
\left.+\int_{\Omega_z^2}\partial_x\left(A_x \partial_xw_1^\prime\left(x,z_\iota\right)\right) \mathrm{d}z+\int_{\Omega_z^2}\partial_x\left(A_x \partial_x\xi_2^\prime\right)\mathrm{d}z\right),
\end{split}
\end{equation}
which is finally equivalent to
\begin{equation}
\begin{split}
%\int_{\Omega_z^2}\partial_x\left(A_x\partial_x\overline{w}\right)=\frac{A_z\left(x,z_\iota\right)}{\varepsilon\left(z_\iota\right)}\partial_z\xi^\prime_2\left(x,z_\iota\right)\\
%+\int_{\Omega_z^2}\partial_x\left(A_x \partial_xw_1^\prime\left(x,z_\iota\right)\right) ~dz+\int_{\Omega_z^2}\partial_x\left(A_x \partial_x\xi_2^\prime\right)~dz
%
\frac{A_z\left(x,z_\iota\right)}{\varepsilon\left(z_\iota\right)}\partial_z\xi^\prime_2\left(x,z_\iota\right)=-\int_{\Omega_z^2}\partial_x\left(A_x \partial_xw\left(x,z_\iota\right)\right) \mathrm{d}z -\int_{\Omega_z^2}\partial_x\left(A_x \partial_x\xi_2^\prime\right)\mathrm{d}z.
\end{split}
\end{equation}

Using this last equation in (\ref{eq:esti1_on_w}), permits to get
\begin{equation}
\begin{split}
|| \sqrt{A_x}\, \partial_x w ||^2_{L^2_1}+ || \sqrt{A_z \over \varepsilon\left(z\right) }\, \partial_z w ||^2_{L^2_1}- \left( \partial_x\left(A_x \partial_x w\left(\cdot ,z_\iota\right)\right) -  \partial_x \left(A_x \partial_x \xi_2^\prime\right)    , w\left(\cdot ,z_\iota\right)\right)_{L^2_2}=0.
\end{split}
\end{equation}
Integrating by parts in the last two terms, yields
\begin{equation}
\begin{split}
|| \sqrt{A_x}\, \partial_x w ||^2_{L^2_1}+ || \sqrt{A_z \over \varepsilon\left(z\right) }\, \partial_z w ||^2_{L^2_1}+|| \sqrt{A_x} \partial_x w\left( \cdot ,z_\iota\right)||^2_{L^2_2} =- \left( A_x \partial_x\xi_2^\prime, \partial_x w\left( \cdot ,z_\iota\right)\right)_{L^2_2}\,,
\end{split}
\end{equation}
which permits to obtain, as ${1 \over \eps(z)} \ge {1 \over \eps_{max}}$ in $\Omega_1$,
\begin{equation}
\begin{split}
\|\partial_xw\|^2_{L^2\left(\Omega_1\right)}+\|\partial_zw\|^2_{L^2\left(\Omega_1\right)}+\|\partial_xw\left(\cdot,z_\iota\right)\|^2_{L^2\left(\Omega_x\right)}\leq c \|\partial_x\xi_2^\prime\|_{L^2\left(\Omega_2\right)}\|\partial_xw\left(\cdot,z_\iota\right)\|_{L^2\left(\Omega_x\right)},
\end{split}
\end{equation}
leading thanks to Theorem \ref{th_xi2} to the $H^1$-estimate of $w$
\begin{equation} \label{AAA}
\begin{split}
\|\partial_xw\|^2_{L^2\left(\Omega_1\right)}+\|\partial_zw\|^2_{L^2\left(\Omega_1\right)}+\|\partial_xw\left(\cdot,z_\iota\right)\|^2_{L^2\left(\Omega_x\right)}\leq c\|\partial_x\xi_2^\prime\|^2_{L^2\left(\Omega_2\right)}\leq c\varepsilon\left(z_\iota\right).
\end{split}
\end{equation}
{\bf 2nd Step}: $H^1$-estimate of $\overline{w}$ and $w_1^\prime$:\\
To estimate from this last inequality the functions $\overline{w}$ and $w_1^\prime$, we shall use the constraint. Indeed, one has
$$
\int_{\Omega_z^1}w(x,z)\left(x,z\right)\mathrm{d}z+L_z^2w \left(x,z_\iota\right)+\int_{\Omega_z^2}\xi_2^\prime\left(x,z\right)\mathrm{d}z= L_z\, \bar{w}(x)\,, \quad \forall x \in \Omega_x\,,
$$
implying
\be \label{BBB}
||\bar{w}||_{L^2(\Omega_x)}  \le c \left( ||w||_{L^2(\Omega_1)} + ||w(\cdot,z_\iota)||_{L^2(\Omega_x)} + ||\xi_2^\prime||_{L^2(\Omega_2)}\right)\leq c\, \sqrt{\varepsilon\left(z_\iota\right)}\,.
\ee
As furthermore $w_1^\prime=w-\bar{w}$, we can conclude the proof in the regular case. In the less regular case, Theorem \ref{th_xi2} as well as formul\ae~(\ref{AAA}), (\ref{BBB}) permit to get the convergences in $H^1$.
\finproof\\
Estimates \eqref{eq:est_on_xi2} and \eqref{ESS} will permit in the following to measure the error done by using the computationally more advantageous (AP/L)-hybrid model instead of a full AP-model for the resolution of \eqref{P}, and to place the interface-position $z_\iota$. 

%Now, we have
%\begin{equation}
%\overline{w}=\frac{1}{L_z}\int_{\Omega_z}w=\frac{1}{L_z}\int_{\Omega_z^1}w+\frac{L_z^2}{L_z}w\left(\cdot,z_\iota\right)
%\end{equation}
%and there exists a positive constant denoted by $C$ too such that
%\begin{equation}
%\|\partial_x\overline{w}\|\leq C\|\partial_xw\|+C\|\partial_xw\left(\cdot,z_\iota\right)\|\leq C\sqrt{\varepsilon\left(z_\iota\right)}.
%\end{equation}
%%
%Moreover, there exists a positive constant $C$ such that
%\begin{equation}
%\|\partial_xw_1^\prime\|\leq\|\partial_xw\|+\|\partial_x\overline{w}\|\leq C\sqrt{\varepsilon\left(z_\iota\right)}
%\end{equation}
%and
%\begin{equation}
%\|\partial_zw_1^\prime\|\=\|\partial_xw\|\leq C\sqrt{\varepsilon\left(z_\iota\right)}.
%\end{equation}
%%
%Poincar\'e inequality let us conclude that $\overline{w}$ and $w_1^\prime$ tends to zero in $H^1$ as $\varepsilon\left(z_\iota\right)$ tends to zero.

%%%%%%%%%%%%%%%%%%%%%%%%%%%%%%%%%%%%%%%%%%%%%%%
\section{Numerical discretization and investigations of the hybrid model} \label{SEC4}
%%%%%%%%%%%%%%%%%%%%%%%%%%%%%%%%%%%%%%%%%%%%%%%

The aim of this section is the discretization of the (AP/L)-coupling model (\ref{weakform}) by means of a finite element method, and the discussion of the obtained numerical results.

%%%%%%%%%%%%%%%%%%%%%%%%%%%%%%%%%%%%%%
\subsection{Finite element discretization}\label{SSEC41}
%%%%%%%%%%%%%%%%%%%%%%%%%%%%%%%%%%%%%%
Let us present here in some details the finite element method we used.

The computational domain $\Omega_x\times\Omega_z$ with $\Omega_x=\left[x_-,x_+\right]$ and $\Omega_z=\left[z_-,z_+\right]$
%, with $L_z=z_+-z_-$,
 is decomposed into a set of rectangular cells $\left[x_i,x_{i+1}\right]\times\left[z_k,z_{k+1}\right]$ where
\begin{eqnarray*}
x_i=x_- + i\Delta x,&&i=0,\dots,N_x+1\,, \qquad \Delta x = \left(x_+-x_-\right) / \left(N_x+1\right)\,,\\
z_k=z_- + k\Delta z,&&k=0,\dots,N_z+1\,, \qquad \Delta z = \left(z_+ -z_-\right) / \left(N_z+1\right)\,.
\end{eqnarray*}
The domain  $\Omega_z$ is decomposed into $\Omega_z^1$ and $\Omega_z^2$, $z_\iota$ with $\iota\in\left\{1,\dots,N_z\right\}$ denoting the interface delimiting these two sub-domains so that $\Omega_z^1$ is composed of the cells $\left[z_k,z_{k+1}\right]$ for $k=\iota,\ldots,N_z$.

The functional spaces $\mathcal{V}$ resp. $\mathcal{V}_1$ are approximated by $\mathbb{Q}_1$ finite element spaces, denoted $\mathcal{V}_{h}$ resp. $\mathcal{V}_{1,h}$. % and defined as
The mean functions and the Lagrangian belonging to $\mathcal{W}$ and $\mathcal{L}$  are approximated by a $\mathbb{P}_1$ finite element  giving rise to the definition of $\mathcal{W}_h$ and $\mathcal{L}_h$. 
In the sequel, $\chi_i(x)$ and $\kappa_k(z)$ denote the standard $\mathbb{P}_1$ hat functions respectively centered in $x_i$ and $z_k$.
% \begin{equation*}
%   \mathcal{W}_h :=  \left\{ \overline{\psi}_h\left(x\right)=\sum_{i=1}^{N_x}\alpha_{i}\chi_i\left(x\right),~x\in\Omega_x\right\} \,,
% \end{equation*}
% and finally, the discrete Lagrangian functional space is defined as
% \begin{equation*}
%  \mathcal{L}_h :=  \left\{ \overline{P}_h\left(x\right)=\sum_{i=0}^{N_x+1}\gamma_{i}\chi_i\left(x\right),~x\in\Omega_x\right\} \,.
% \end{equation*}

The weak formulation of the hybrid (AP/L)-model (\ref{weakform}) is approximated thanks to a three points Gauss quadrature formula, yielding the following linear system
\begin{equation}\label{eq:system:hybrid}\left(AP/L\right)_h~~~~~~~~~~
\left(\begin{array}{ccc}
A_{xa} & \frac{1}{L_z}\left(C_{a1}+C_{a2}^\iota\right)  & 0 \\
C_{f1} & A_{xf1}+A_{z1} & B_{l1} \\
%C_{f2} & -D & A_{xf2}& 0  \\
0 & B_{c1}+B_{c2}^\iota & 0
\end{array}\right)\left(\begin{array}{c}
\alpha  \\
\beta  \\
\gamma %\\
%\delta
\end{array}\right)=\left(\begin{array}{c}
F_{\overline{u}} \\
F_{u_1^\prime} \\
%F_{u_2^\prime} \\
0
\end{array}\right),
\end{equation}
where $A_{\star}$ (resp. $B_{\star}$, $C_{\star}$) is the matrix associated to the bilinear form $a_{\star}$ (resp. $b_{\star}$, $c_{\star}$) introduced in \eqref{eq:bilinear:forms}.
Note that the matrices  $C_{a2},~B_{c2}\in\mathbb{R}^{N_x\times N_x}$ associated to $a_{a2},~b_{c2}$ involve only the nodes located on the interface, such that they are expanded to $C_{a2}^\iota,~B_{c2}^\iota$ by adding zero elements in order to conform with the size of $C_{a1},~B_{c1}\in\mathbb{R}^{N_x\times N_x\left(N_z+2-\iota\right)}$  to define the system \eqref{eq:system:hybrid}. The right-hand side definition is specified here
\begin{eqnarray*}
\left(F_{\overline{u}}\right)_i&:=&\left(\overline{f},\chi_i\right)_{L^2\left(\Omega_x\right)}+\frac{1}{L_z}\left(g_+-g_-,\chi_i\right)_{L^2\left(\Omega_x\right)},~ \forall~i=1,\dots,N_x,\\
\left(F_{u_1^\prime}\right)_{ik}&:=&\left(f,\chi_i\kappa_k\right)_{L^2\left(\Omega_1\right)}+\left(g_+,\chi_i\kappa_k\left(z_+\right)\right)_{L^2\left(\Omega_x\right)},~ \forall~i=1,\dots,N_x,~k=\iota,\dots,N_z+1,%\\
%\left(F_{u_2^\prime}\right)_i&:=&\left(f,\chi_i\right)_{|\Omega_x\times\Omega_z^2}-\left(g_-,\chi_i\right),
\end{eqnarray*}
and the unknowns of the system are three vectors $\alpha\in\mathbb{R}^{N_x}$, $\beta\in\mathbb{R}^{N_x\left(N_z+2-\iota\right)}$ and $\gamma\in\mathbb{R}^{N_x+2}$
\begin{equation*}
\begin{array}{lclclcl}
\alpha&:=&\left(
\alpha_{1}
\cdots
\alpha_{N_x}
\right)^T, &&
\beta&:=&\left(\beta_{1\iota}
\cdots
\beta_{1N_z+1}
\cdots
\beta_{N_x\iota}
\cdots
\beta_{N_xN_z+1}
\right)^T, \\
\gamma&:=&\left(
\gamma_{0}
\cdots
\gamma_{N_x+1}
\right)^T.
\end{array}
\end{equation*}

As a comparison, the linear system corresponding to the AP-model (\ref{weakformAP}) writes
\begin{equation}\left(AP\right)_h~~~~~~~~~~
\left(\begin{array}{ccc}
A_{xa} & \frac{1}{L_z}C_a  & 0 \\
C_f & A_{xf}+A_z & B_l \\
0 & B_c & 0
\end{array}\right)\left(\begin{array}{c}
\alpha  \\
\tilde{\beta}  \\
\gamma 
\end{array}\right)=\left(\begin{array}{c}
F_{\overline{u}} \\
F_{u^\prime} \\
0
\end{array}\right),\label{eq:def:Matrix:AP}
\end{equation}
where $A_z$ (resp. $A_{xf}$, $B_l$, $B_c$, $C_a$ and $C_f$) is associated to the bilinear form similar to $a_{z1}$ (resp. $a_{xf1}$, $b_{l1}$, $b_{c1}$, $c_{a1}$ and $c_{f1}$) but integrated on $\Omega_z$ instead of $\Omega_z^1$. The right-hand side $F_{u^\prime}$ is defined for $i=1,\dots,N_x,~k=0,\dots,N_z+1$, by
\begin{eqnarray*}
\left(F_{u^\prime}\right)_{ik}&:=&\left(f,\chi_i\kappa_k\right)_{L^2\left(\Omega\right)}+\left(g_+,\chi_i\kappa_k\left(z_+\right)\right)_{L^2\left(\Omega_x\right)}-\left(g_,\chi_i\kappa_k\left(z_-\right)\right)_{L^2\left(\Omega_x\right)},%\\
%\left(F_{u_2^\prime}\right)_i&:=&\left(f,\chi_i\right)_{|\Omega_x\times\Omega_z^2}-\left(g_-,\chi_i\right),
\end{eqnarray*}
and the unknown $\tilde{\beta}\in\mathbb{R}^{N_x\left(N_z+2\right)}$ writes
\begin{equation*}
\tilde{\beta}:=\left(\beta_{10}
\cdots
\beta_{1N_z+1}
\cdots
\beta_{N_x0}
\cdots
\beta_{N_xN_z+1}
\right)^T. 
\end{equation*}

The gain of the hybrid (AP/L)-model as compared to the fully AP-model results from the size reduction of the fluctuation unknown. Indeed, $\tilde{\beta}$ is a $N_x\left(N_z+2\right)$-vector, while $\beta$ is a $N_x\left(N_z+2-\iota\right)$-vector.

As for the P-model (\ref{P_var}), the corresponding linear system 
%is constituted of $N_x\left(N_z+2\right)$ equations and 
is of the form
\begin{equation}\left(P\right)_h~~~~~~~~~~
%\left(\begin{array}{c}
(A_{xf}+A_z) 
%\end{array}\right)\left.\begin{array}{c}
\delta 
%\end{array}\right.
=%\left.\begin{array}{c}
F\,,
%\end{array}\right.,
\label{eq:def:system:P}
\end{equation}
where the second member is defined for $i=1,\dots,N_x,~k=0,\dots,N_z+1$, by
\begin{eqnarray*}
\left(F\right)_{ik}&:=&\left(f,\chi_i\kappa_k\right)_{L^2\left(\Omega\right)}+\left(g_+,\chi_i\kappa_k\left(z_+\right)\right)_{L^2\left(\Omega_x\right)}-\left(g_,\chi_i\kappa_k\left(z_-\right)\right)_{L^2\left(\Omega_x\right)},
\end{eqnarray*}
and the unknown $\delta\in\mathbb{R}^{N_x\left(N_z+2\right)}$ writes
\begin{equation*}
\delta:=\left(\delta_{10}
\cdots
\delta_{1N_z+1}
\cdots
\delta_{N_x0}
\cdots
\delta_{N_xN_z+1}
\right)^T. 
\end{equation*}

%\begin{bfseries}
% \noindent Ecrire la systeme lineaire pour le schema AP classique.\\
%  Faire une phrase pour la matrice du probleme P ?
%  Expliquer que le gain de la methode hybride le gain vient du fait que 
% l'on a reduit la taille du vecteur $\beta$ 
%\end{bfseries}
Note that the linear system providing the fluctuation in \eqref{eq:def:Matrix:AP} is nothing but the P-model matrix \eqref{eq:def:system:P} augmented with the two sub-matrices $B_l$ and $B_c$ discretizing the zero mean value constraint.

%\begin{bfseries}
% \noindent Ajouter une reference pour MUMPS.
%\end{bfseries}
%%%%%%%%%%%%%%%%%%%%%%%%%%%%%%%%%%%%%%
\subsection{Numerical investigations}\label{SSEC42}
%%%%%%%%%%%%%%%%%%%%%%%%%%%%%%%%%%%%%%

\subsubsection{Test case setup}\label{sec:setup}

The efficiency of the hybrid model is illustrated in this section. For this an anisotropy ratio with large variations within the computational domain is constructed thanks to the following definition 
\begin{equation}
\label{eq:defeps}
\varepsilon\left(z\right)=\frac{1}{2}\Big(\varepsilon_{max}\big(1+\tanh \left(rz\right)\big)+\varepsilon_{min}\big(1-\tanh \left(rz\right)\big)\Big),
\end{equation}
where $\varepsilon_{min}>0$ and $\varepsilon_{max}>0$ %\in\mathbb{R}^{+\star}>0 
define the range of values covered by $\varepsilon$, $r \in \mathbb{R}$ defining the rate of change. In the sequel this parameter value is set to $r=30$. The different heterogeneous anisotropy ratios used for the numerical investigations are represented on Figure~\ref{fig:epsilon} as a function of $z$. 
\begin{figure}[!ht]
\begin{center}
\includegraphics[angle=-90,width=0.49\textwidth]{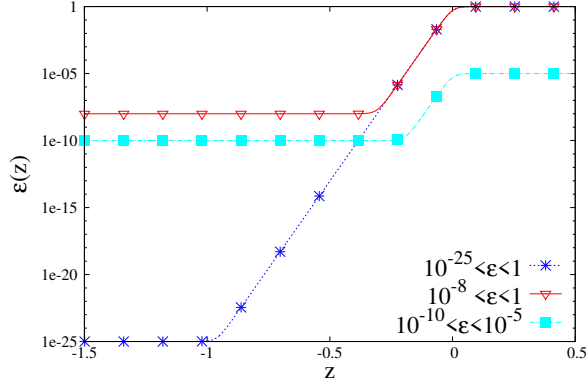}
\caption{Anisotropy ratio $\varepsilon$ as a function of $z$, for different values of $\varepsilon_{min}$ and $\varepsilon_{max}$ with $r=30$.}
\label{fig:epsilon}
\end{center}
\end{figure}

The diffusion coefficients are defined as 
\begin{equation}\label{eq:def:coeffs}
A_x\left(x,z\right)=c_1+xz^2,~~~A_z\left(x,z\right)=c_2+xz,%\\
%A_o\left(x,z\right)&=&c_1+xz^2,
\end{equation}
with $c_1\in\mathbb{R}$, $c_2\in\mathbb{R}$ two constants chosen to meet the requirements $A_x\left(x,z\right)\geq m_x>0$ and $A_z\left(x,z\right)\geq m_z>0$ for $\left(x,z\right) \in \Omega$.

With these definitions an analytic setup is manufactured thanks to the exact solution 
\begin{equation}\label{eq:sol:exact}
u_e\left(x,z\right)=\sin\left(\frac{2\pi}{L_x}x\right)\left(1+\varepsilon\left(z\right) \sin\left(\frac{2\pi}{L_z}z\right)\right)\,,
\end{equation}
by inserting this expression into (\ref{eq:SP})-(\ref{eq:SP:tensor}) for the computations of the second member of the problem $f$ as well as $g_+$ and $g_-$.
A second setup is investigated with the following definitions for the diffusion coefficients
\begin{equation}\label{eq:def:coeffs:b}
A_x\left(x,z\right)=1+\cos(c_1+xz),~~~A_z\left(x,z\right)=1 + \sin^2(c_2+xz)\,,
\end{equation}
and the exact solution
\begin{equation}\label{eq:sol:exact:b}
  u_e\left(x,z\right)=\sin\left(\frac{2\pi}{L_x}x\right)\left(1+ \sin\left(\frac{2\pi}{L_z} \varepsilon\left(z\right)z\right)\right)\,.
\end{equation}
 The simulations are performed with $c_1=c_2=L_z$, on the domain $\Omega_a=\Omega_x\times\Omega_z=\left[0,1\right]\times\left[-1,1\right]$ for homogeneous $\varepsilon$ and with $\Omega_b=\Omega_x\times\Omega_z=\left[0,1\right]\times\left[-1.5,0.5\right]$ for $\varepsilon$ defined by (\ref{eq:defeps}). The results obtained for both setups \eqref{eq:def:coeffs}-\eqref{eq:sol:exact} and \eqref{eq:def:coeffs:b}-\eqref{eq:sol:exact:b} being very similar, only the computations performed thanks to the first definitions are reported in the sequel.

 The linear systems stemming from the discretization of the singular perturbation problem (\ref{P_var}), the AP-model (\ref{weakformAP}) and the hybrid (AP/L)-model (\ref{weakform}) are solved using the same solver. For the computations carried out in this paper, the sparse direct solver MUMPS \cite{MUMPS} is used.
 %%%%%%%%%%%%%%%%%%%%%%%%%%%%%%%%%%%%%%
\subsubsection{Choice of the interface}\label{sec:interface:choice}
%%%%%%%%%%%%%%%%%%%%%%%%%%%%%%%%%%%%%%

%An important difficulty of our coupling is the choice of $z_\iota$. The position of the interface has indeed a big influence on the performance of the scheme, in terms of precision and computational time.

%\textit{\textbf{\`A compl\'eter avec l'analyse math\'ematique.}}

The choice of the interface location $z_\iota$ is of great importance because it influences the performance of the scheme, in terms of precision and computational time. 

This question is often pointed out in hybrid models based on domain decomposition strategies.  
As explained in the introduction and the beginning of section \ref{SEC3}, the motivation of our (AP/L)-coupling strategy is that it is always possible to situate the interface, which is not always the case in a (P/L)-coupling method.
The position of the interface has to be chosen on one hand in order to preserve the precision of the numerical scheme and on the other hand to achieve a gain in the computational time ({\it i.e.} largest possible $L$-sub-domain). This means that the approximation error between the AP-solution and the hybrid (AP/L)-solution has to be of the order of the precision of the numerical scheme.

%In our context, we can for example notice that the coupling of the L-model directly to the P-one is not a good strategy, because it is not always possible to find a region in the computational domain where both models (P) as well as (L) provide accurate approximations of the solution. Degond \textit{et al.} show indeed in \cite{DDLNN} that the domain of accuracy of the L-model does not always intersect the one of the P-model (see Figure~3.1 of \cite{DDLNN}). This motivated our (AP/L)-coupling, since the AP-model is usable in the whole domain. As soon as $\varepsilon$ is small in one part of the domain, it is always possible to locate the interface in a part of the domain, where $\varepsilon$ is small enough such that the precision of the numerical scheme is preserved. 

The analytic investigations carried out in section~\ref{SEC3} (see Theorem~\ref{TheTheorem}) demonstrate that the error produced by the use of the limit problem in one sub-domain is bounded by the value of the asymptotic parameter at the interface, {\it i.e.} $\sqrt{\varepsilon(z_{\iota})}$. If we denote by $u$ the exact solution of the problem, by $v$ the exact solution of the hybrid model and   finally by $v_h$ its numerical approximation, the following inequality $||u - v_h ||_{H^1} \leq ||u- v ||_{H^1} + || v-v_h||_{H^1}$ holds, leading to
  \begin{equation}\label{eq:error:estimate}
    ||u - v_h ||_{H^1} \leq c\sqrt{\varepsilon(z_{\iota})} + C h^m\,,
  \end{equation}
$C >0$ denoting a constant, $h$ the typical mesh size and $m$ the approximation order of the numerical methods defined as the convergence rate of the error in the $H^1$-norm.
 In order to prevent any deterioration of the numerical method precision, the interface should hence be located in a region of the computational domain where the following condition is met $\sqrt{\varepsilon(z_{\iota})}\sim h^m$.

% \begin{figure}[!ht]
%   \centering
%   \includegraphics[angle=-90,width=0.6\textwidth]{FiguresReview/errorH1etL2tc8Nqp3relDa.eps}
%   \caption{Relative error evolution as a function of the asymptotic parameter value at the interface. Error norm computed between a reference solution carried out thanks to the AP-scheme and the approximation computed by the (AP/L)-scheme.}
%   \label{fig:interface}
% \end{figure}
To verify numerically the statements of Theorem~\ref{TheTheorem} and more particularly the error estimate provided by the equation~\eqref{eq:error:estimate}, a series of computations have been run.
% The first outputs displayed on the Figure~\ref{fig:interface} relate the error $||u - v_h ||_{H^1}$ with the interface value $\varepsilon(z_{\iota})$.} These computations are carried out on a $512\times 512$ mesh, with an anisotropy ratio defined by equation \eqref{eq:defeps} where $\varepsilon_\text{min}= 10^{-25}$ and $\varepsilon_\text{max}=1$. The right hand side of problem \eqref{P} is defined by
% \begin{equation}
%   f(x,z)=1+\exp( - ( 16\, (x-1/2)^2 + 8\,z^2 ) )\,,% \quad g_{\pm}=0 \,,
% \end{equation}
% and the diffusion coefficients $A_x$ and $A_z$ are given by \eqref{eq:def:coeffs}. For this set up it is unfortunately not possible to explicit the exact solution of the problem, thus a reference solution is computed thanks to the full AP-scheme \eqref{weakformAP} used in the whole computational domain on a mesh of $512\times 512$ cells. What we observe from  Figure~\ref{fig:interface} is a linear decrease of the hybrid method approximation error (between the numerical and reference solutions) with respect to $\varepsilon(z_{\iota})$.
%
% without any success in reproducing . The test cases we chose are apparently so particular, that the  $\sqrt{\varepsilon(z_{\iota})}$-estimate is too pessimistic. However, to be sure not to encounter less gentle cases, the interface has to be located such that the condition $\sqrt{\varepsilon(z_{\iota})}\sim h^m$ is met.
The evolution of the error measured between the exact solution (manufactured thanks to the set-up \eqref{eq:def:coeffs} and \eqref{eq:sol:exact}) and the approximation carried out thanks to the hybrid method is displayed on Figure~\ref{fig:optimal} (left plots). This Figure displays also  the value of the optimal asymptotic parameter at the interface denoted  $\varepsilon(z_{\iota}^\star)$ as a function of the mesh size. The optimal value of $\varepsilon(z_{\iota})$ is the largest value preserving the accuracy of the numerical method.  As reported on  Figure~\ref{fig:optimal}, the approximation error is decreasing with $\varepsilon(z_{\iota})$ until a plateau is reached. 
\begin{figure}[!ht]
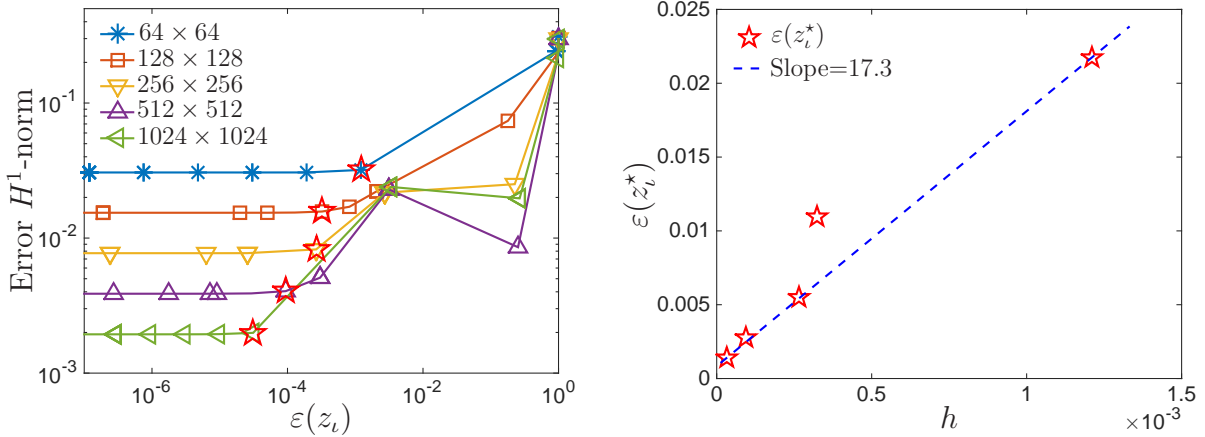

  \centering
  {\psfrag{1}[l][l][.8]{$64\times 64$}
  \psfrag{2}[l][l][.8]{$128\times 128$}
  \psfrag{3}[l][l][.8]{$256\times 256$}
  \psfrag{4}[l][l][.8]{$512\times 512$}
  \psfrag{5}[l][l][.8]{$1024\times 1024$}
  \psfrag{X}[][][1.]{$\varepsilon(z_\iota)$}
  \psfrag{Y}[][][1.]{Error $H^1$-norm}
  \includegraphics[width=0.48\textwidth]{Error.eps}\hspace*{0.04\textwidth}}
  \psfrag{Data}[l][l][.8]{${\varepsilon(z_\iota^\star)}$}
  \psfrag{Slope=17.27}[l][l][.8]{Slope=17.3}
  \psfrag{Y}[][][1.]{$\varepsilon(z_\iota^\star)$}
  \psfrag{X}[][][1.]{$h$}
  \includegraphics[width=0.48\textwidth]{EpsOpt.eps}
  \caption{Error $H^1$-norm between the exact solution and the approximation computed by the (AP/L)-scheme on different meshes (left), the optimal value of $\varepsilon(z_\iota)$ denoted $\varepsilon(z_\iota^\star)$ is identified by a pentagram marker; Optimal value of ${\varepsilon(z_\iota)}$  as a function of the mesh size $h$ (right).}
  \label{fig:optimal}
\end{figure}
 Above this $\varepsilon(z_\iota^\star)$-value, the error is dominated by the modeling error, due to the use of the limit model for large values of the asymptotic parameter. Below this value, the error does not depend on $\varepsilon(z_{\iota})$ anymore and reduces to the numerical approximation error proportional to $h^m$. This optimal value ensures that the limit model is used on the largest sub-domain as possible without introducing any discrepancy in the numerical method precision. The second plot displays $\varepsilon(z_{\iota}^\star)$ as a function of the mesh size. A linear decrease of the ${\varepsilon(z_{\iota}^\star)}$ values with the mesh size is observed, analogously to the decrease of the error $H^1$-norm with the ${\varepsilon(z_{\iota})}$-values, while this decrease is expected, as proved in Theorem~\ref{TheTheorem} and stated in the equation~\eqref{eq:error:estimate} to be proportional to $\sqrt{\varepsilon(z_{\iota})}$.
 This test case has been reproduced with different set-ups obtaining the same decrease rate of the approximation error. However, the $\sqrt{\varepsilon(z_{\iota})}$-estimate of  Theorem~\ref{TheTheorem} is the only estimate we can guarantee such that the condition $\sqrt{\varepsilon(z_{\iota})}\sim h^m$ should be used to locate the interface.

 This hybrid method has been developed to accelerate the simulation of time dependent systems described by a set of equations, in which the anisotropic elliptic equation investigated here provides one of the unknowns. This system is assumed to incorporate two different time scales. The first one is related to the evolution of the quantities advanced by the set of equations. The second one, much slower, characterizes the variations of the anisotropy ratio. This is the typical framework for the simulation of ionospheric plasma disturbances already mentioned in the introduction. Indeed the dynamics of the plasma perturbations is fast compared to that of the neutral particles properties that are responsible for the anisotropy variations. In this framework, the interface can be precisely located at the initial time, thanks to the use of the AP-scheme providing a reference solution. Note that the implementation of the AP-scheme is readily obtained from that of the hybrid model. Then the computation of a series of time steps can be accelerated using the hybrid model for the anisotropic equation resolution. When the anisotropy ratio value undergoes significant variations, the procedure can be repeated to determine a new location for the interface. The hybrid method is thus all the more efficient than the time scales separation is large.
\begin{remark}
The hybrid model requires two sub-domains verifying: $\Omega_z^1\neq \emptyset$ and $\Omega_z^2\neq \emptyset$. This means that for $\varepsilon \sim \mathcal{O}(1)$ uniformly in the domain, the method introduced here cannot provide accurate computations since the limit problem will be used in one part of the domain where the asymptotic parameter is large. In such a situation the P-model or the AP-formulation should be preferred.

%It is necessary to have two subdomains: $\Omega_z^1\neq \emptyset$ and $\Omega_z^2\neq \emptyset$. If $\varepsilon\left(z\right)\leq h^2~\forall z\in\Omega_z$, we choose $\iota=N_z$. If $\varepsilon\left(z\right)> h^2~\forall z\in\Omega_z$, we choose $\iota=1$.
\label{rk:nonemptyset}
\end{remark}

%%%%%%%%%%%%%%%%%%%%%%%%%%%%%%%%%%%%%%
\subsubsection{Accuracy of the numerical method}
%%%%%%%%%%%%%%%%%%%%%%%%%%%%%%%%%%%%%%

The precision of the hybrid numerical method is first examined. 
Primarily, we represent the solution computed by the (AP/L)-scheme when $\varepsilon$ is given by (\ref{eq:defeps}) with $\varepsilon_{min}=10^{-8}$, $\varepsilon_{max}=1$, $N_x=N_z=64$ on Figure~\ref{fig:sol3dtanh}. We also give the difference between this computed solution and the exact one on Figure~\ref{fig:err_sol3dtanh}.
%\begin{figure}[!ht]
%\begin{center}
%\subfigcapskip -0.5em
%\psfrag{sqrt}[][][1.]{h}
%\subfigure[\label{fig:sol3dcst}]{\includegraphics[angle=-90,width=0.49\textwidth]{FiguresReview/u_APL_tc3_epsmin1e-8_Nx64_Nz64_Nqp3.eps}}
%\subfigure[\label{fig:err_sol3dcst}]{\includegraphics[angle=-90,width=0.49\textwidth]{FiguresReview/error_APL_tc3_epsmin1e-8_Nx64_Nz64_Nqp3.eps}}
%\caption{Constant $\varepsilon=10^{-8}$ case: (a) solution of the (AP/L)-scheme, (b) difference with the exact solution.}
%\label{fig:slope}
%\end{center}
%\end{figure}

\begin{figure}[!ht]
\begin{center}
\subfigcapskip -0.5em
\psfrag{sqrt}[][][1.]{h}
\subfigure[\label{fig:sol3dtanh}]{\includegraphics[angle=-90,width=0.49\textwidth]{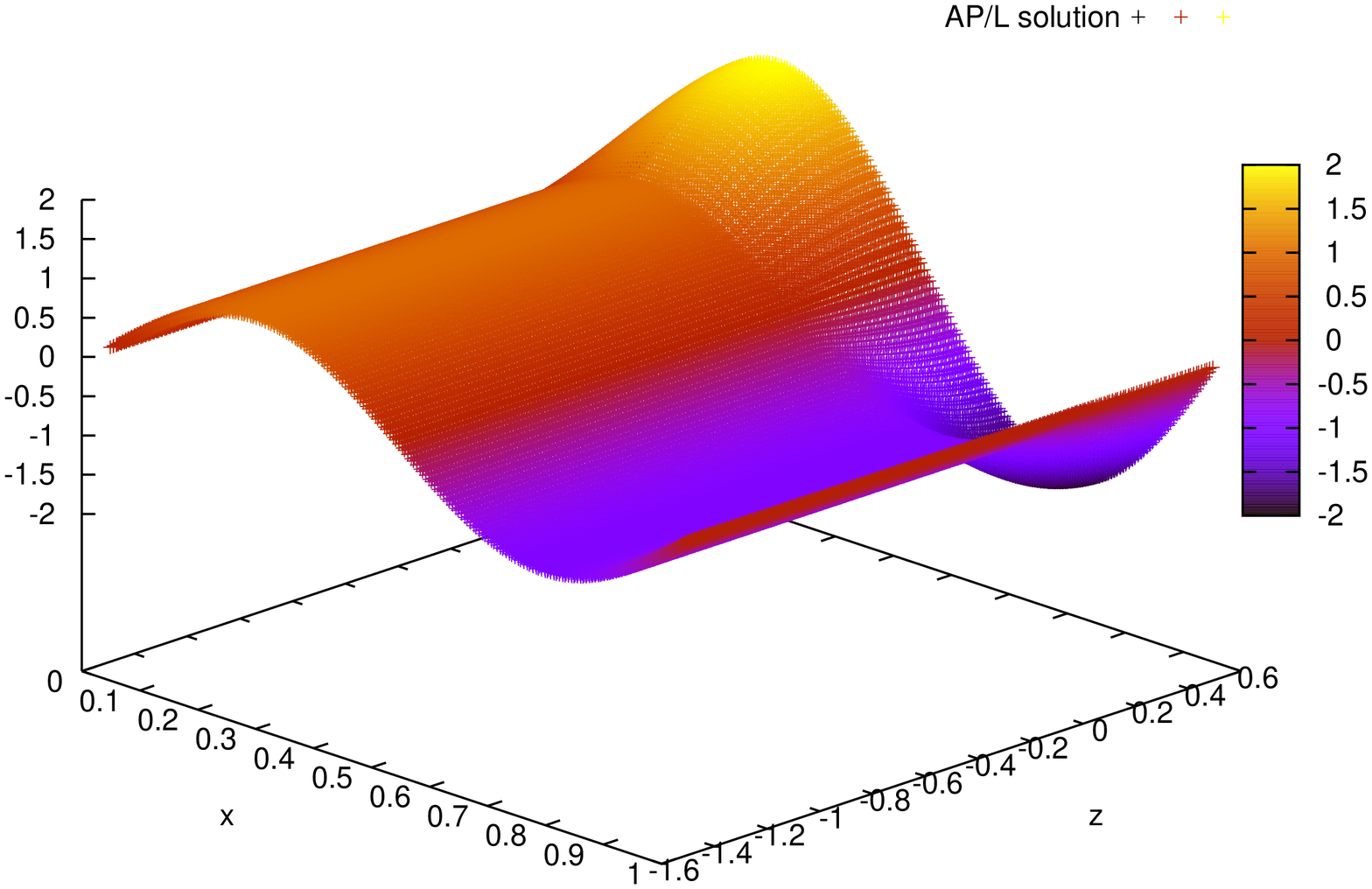}}
\subfigure[\label{fig:err_sol3dtanh}]{\includegraphics[angle=-90,width=0.49\textwidth]{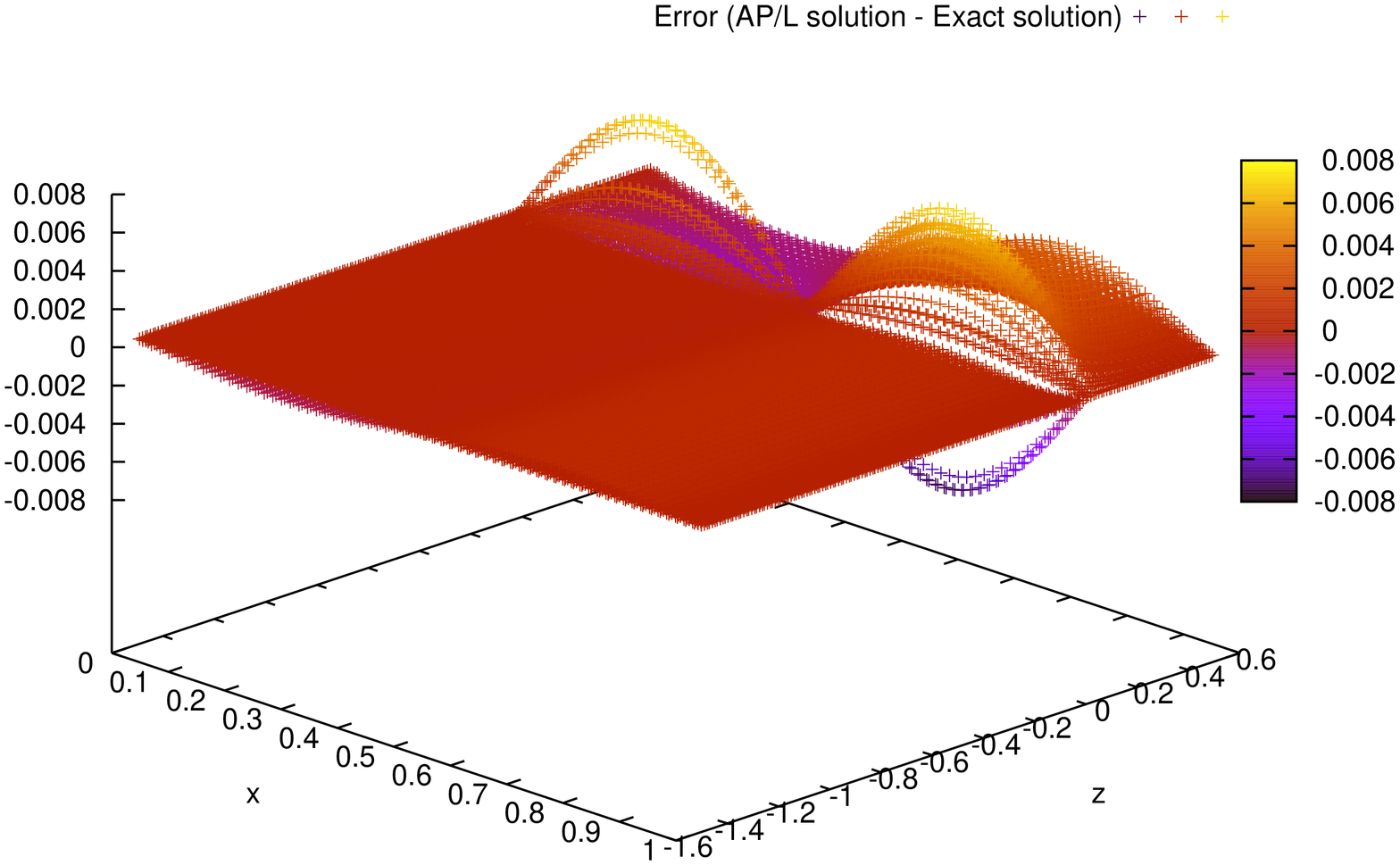}}
\caption{$\varepsilon$ given by (\ref{eq:defeps}) case, with $\varepsilon_{min}=10^{-8}$, $\varepsilon_{max}=1$: (a) solution of the (AP/L)-scheme, (b) difference with the exact solution.}
\label{fig:slope}
\end{center}
\end{figure}

The relative H$^1$-error between the exact solution \eqref{eq:sol:exact} and its numerical approximation via the (AP/L)-model is plotted on Figure~\ref{fig:slopecst} as a function of the mesh size $h:=(\Delta x\Delta z)^\frac{1}{2}$. Note that the computations are performed with the same number of cells in each direction, which amounts, accordingly to the setup precised in section \ref{sec:setup}, to $\Delta z = 2 \Delta x$. These errors are computed first, for constant $\varepsilon$, with $\varepsilon$ ranging from $10^{-25}$ to $10^{-8}$  (Figure~\ref{fig:slopecst}). Second, the same computations are performed and plotted on Figure~\ref{fig:slopevar} for heterogeneous anisotropy ratios $\varepsilon$ defined by \eqref{eq:defeps} with $\varepsilon_{max}$ set to $1$ and $\varepsilon_{min}$ ranging from $10^{-25}$ to $10^{-8}$  as displayed on Figure~\ref{fig:epsilon} (curves with triangle and star markers). The reference solution is manufactured thanks to the setup (\ref{eq:def:coeffs})-(\ref{eq:sol:exact}).
\begin{figure}[!ht]
\begin{center}
\subfigcapskip -0.5em
\psfrag{sqrt}[][][1.]{h}
\subfigure[\label{fig:slopecst}]{\includegraphics[angle=-90,width=0.49\textwidth]{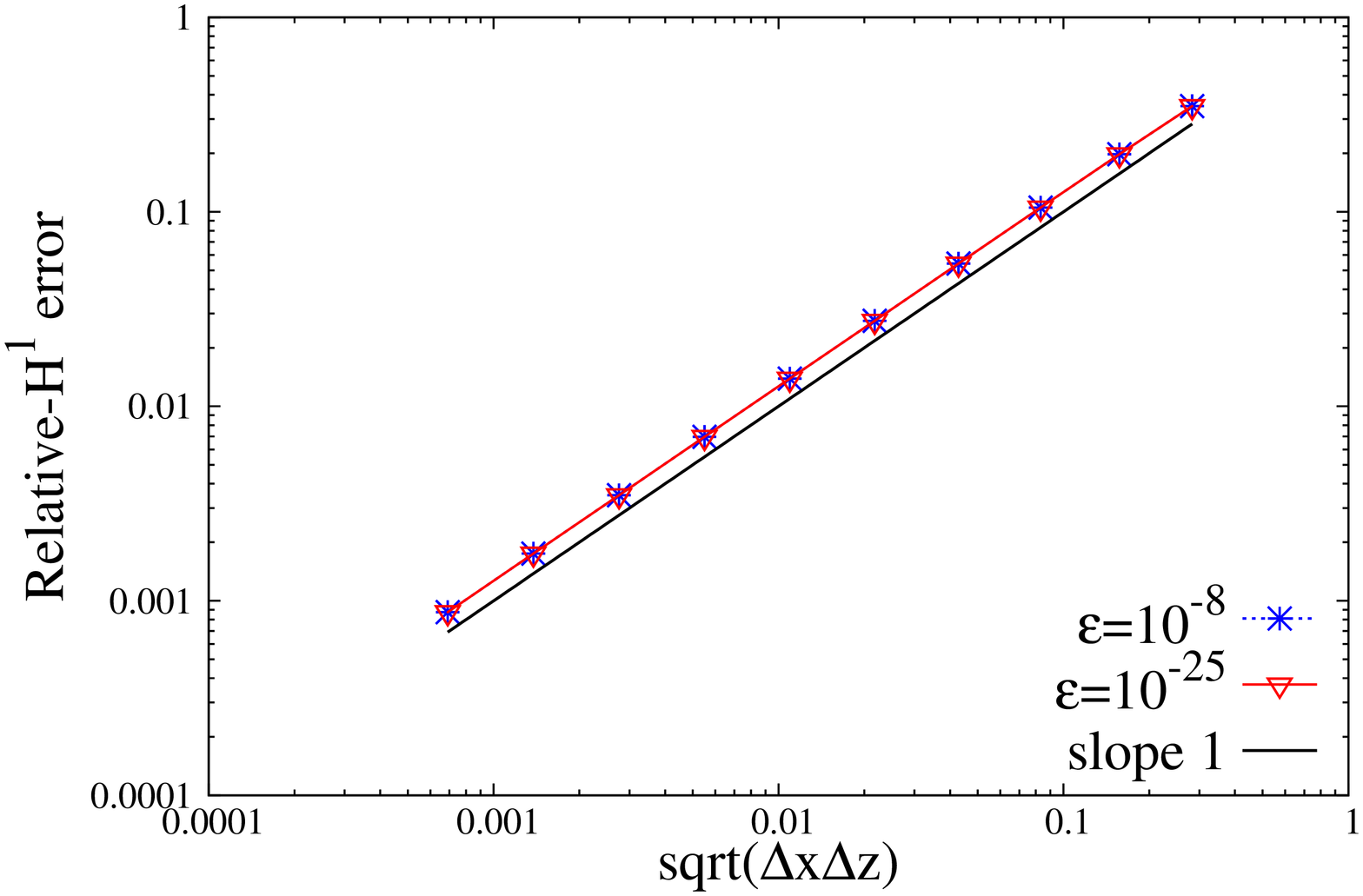}}
\subfigure[\label{fig:slopevar}]{\includegraphics[angle=-90,width=0.49\textwidth]{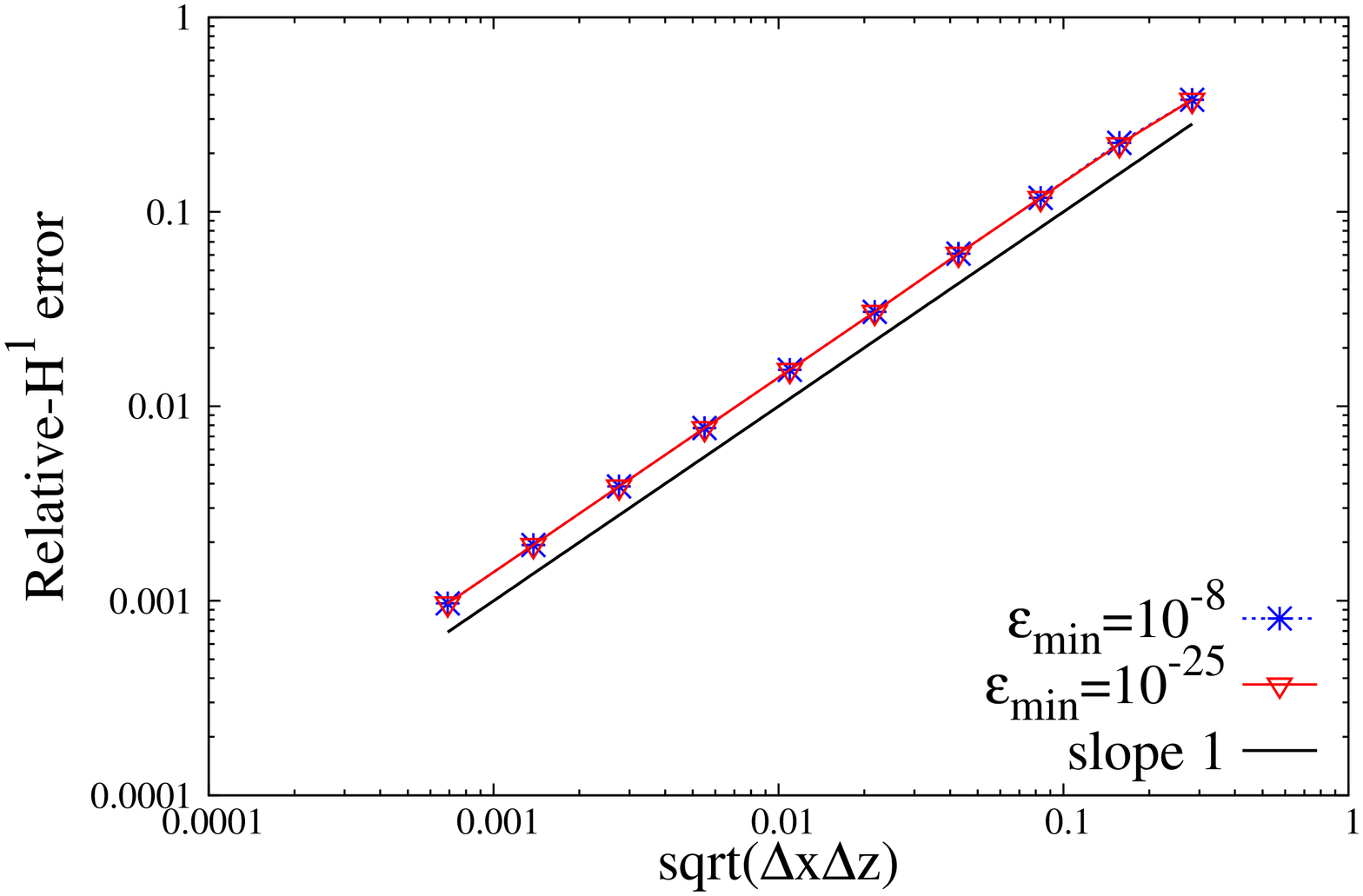}}
\caption{H$^1$ relative error between the exact solution and its numerical approximation as a function of the mesh size $h:=(\Delta x\Delta z)^\frac{1}{2}$: (a) accuracy for constant anisotropy ratios $\varepsilon= 10^{-8}$ and $\varepsilon=10^{-25}$, (b) accuracy for heterogeneous anisotropy ratios defined by \eqref{eq:defeps} with $\varepsilon_\text{max}= 1$, and $\varepsilon_\text{min}= 10^{-8}$ or $\varepsilon_\text{min}=10^{-25}$.}
\label{fig:slope}
\end{center}
\end{figure}
%
%The L$^2$ relative errors computed with heterogeneous ratios are plotted on Figure~\ref{fig:slopevar}.

The hybrid method is observed to be first order accurate for both sets of computations as soon as the interface is located in a region where the value of $\varepsilon$ is small enough compared to the precision of the spatial discretization.

\subsubsection{Asymptotic preserving property of the hybrid model}
\label{sec:num:AP}

The Asymptotic-Preserving property of the hybrid model is now investigated. The error between the solution and its numerical approximations computed thanks to the hybrid (AP/L)-model (\ref{weakform}), the standard AP-scheme (\ref{weakformAP}) and the discretized singular perturbation problem (\ref{P_var}), as well as the condition number of the corresponding linear systems, are plotted on Figures~\ref{fig:ap64} and \ref{fig:ap1024}. These computations are carried out on meshes with $64\times64$ and $1024\times 1024$ cells with the heterogeneous anisotropy ratio defined by \eqref{eq:defeps}. The condition numbers reported on these figures are estimated by the MUMPS solver. Concerning the condition number, remark that the relative error of the linear system solution, carried out by the solver, is upper bounded by two parameters denoted by CN$_1$ and CN$_2$ (see subsection~2.3 of \cite{MUMPS}) and plotted on Figures~\ref{fig:ap64:b} and \ref{fig:ap1024:b}. %The computations being carried out using the double precision arithmetic (with roughly 15 significant digits), the accuracy of the linear system solution is guaranteed for values of these two parameters small compared to $10^{15}$. Note that the CN$_1$ parameter being equal to one for all the computations performed with both the (AP) and the hybrid (AP/L)-schemes, the plots of this parameter for these two methods are omitted on Figures~\ref{fig:ap64:b} and \ref{fig:ap1024:b}.
\begin{figure}[!ht]
\begin{center}
\subfigcapskip -0.5em
\subfigure[\label{fig:ap64:a}]{\includegraphics[angle=-90,width=0.49\textwidth]{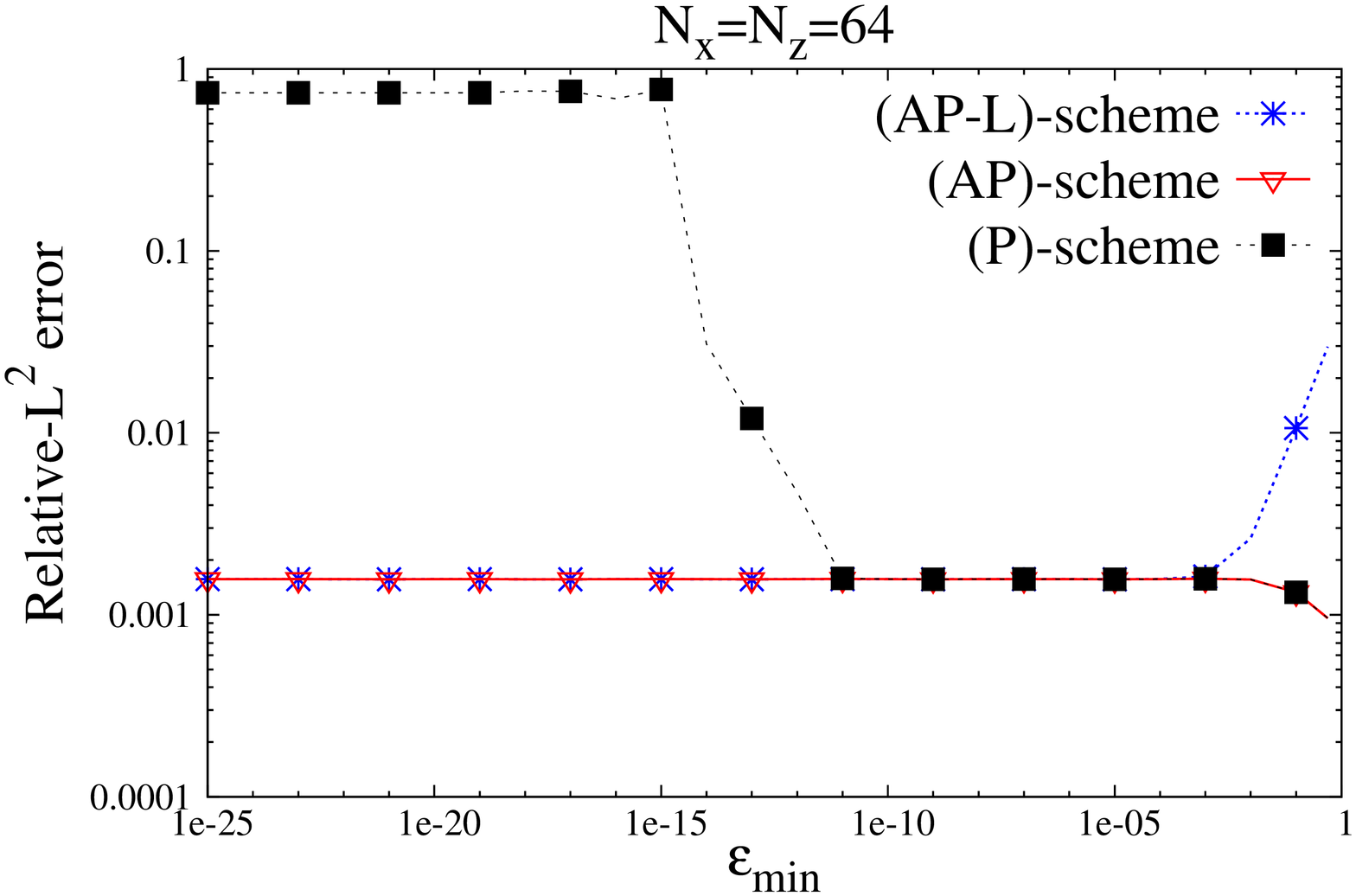}}
\subfigure[\label{fig:ap64:b}]{\includegraphics[angle=-90,width=0.49\textwidth]{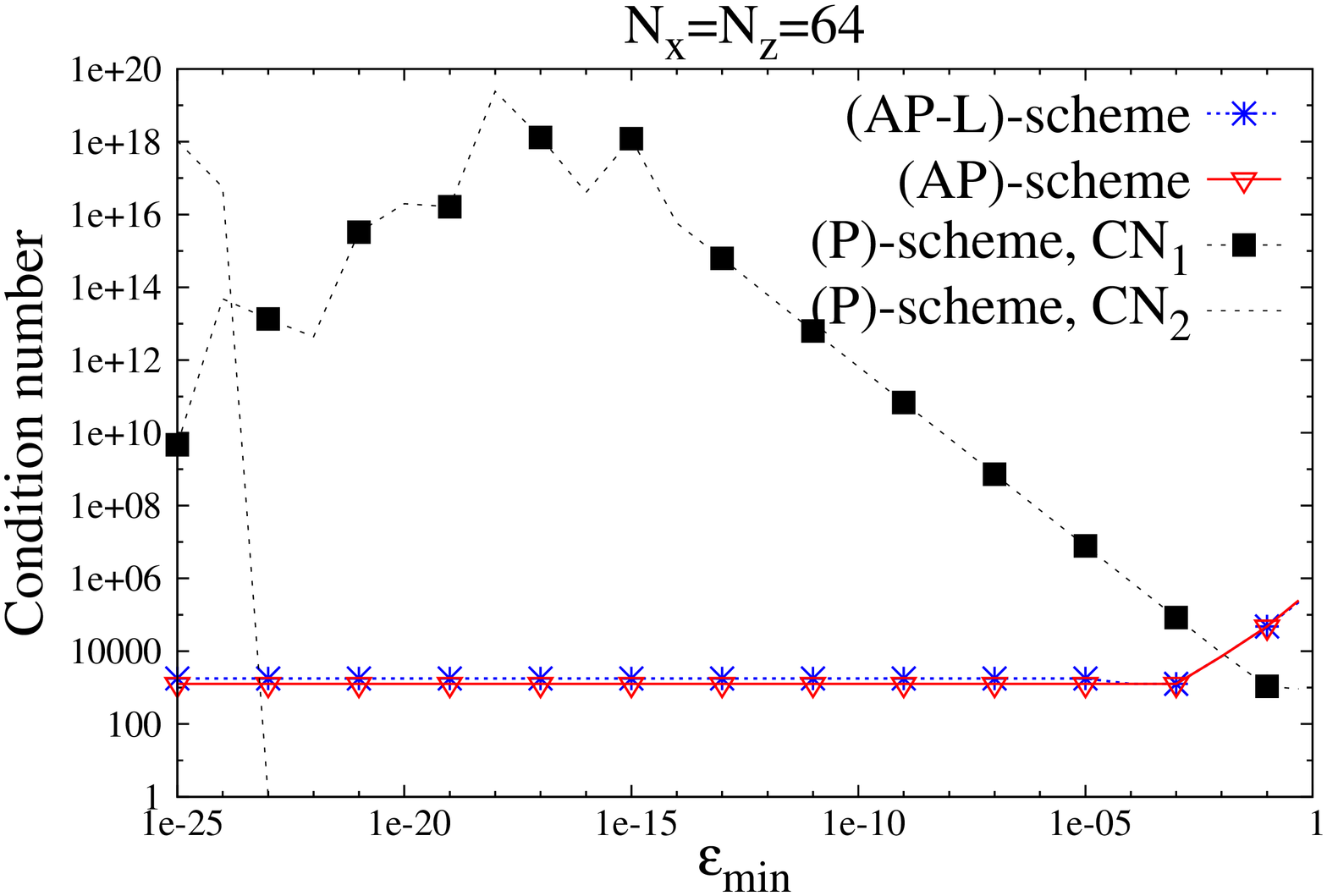}}
\caption{Comparison of the hybrid (AP/L)-model, the AP-scheme and the singular perturbation P-model accuracy for heterogeneous anisotropy ratios as defined by \eqref{eq:defeps} on a $64\times64$ mesh: (a) relative L$^2$-error between the exact solution and its numerical approximations and (b) condition number of the linear systems (estimated by MUMPS) as a functions of $\varepsilon_{min}$ (with $\varepsilon_{max}=1$).}
\label{fig:ap64}
\end{center}
\end{figure}

\begin{figure}[!ht]
\begin{center}
\subfigcapskip -0.5em
\subfigure[\label{fig:ap1024:a}]{\includegraphics[angle=-90,width=0.49\textwidth]{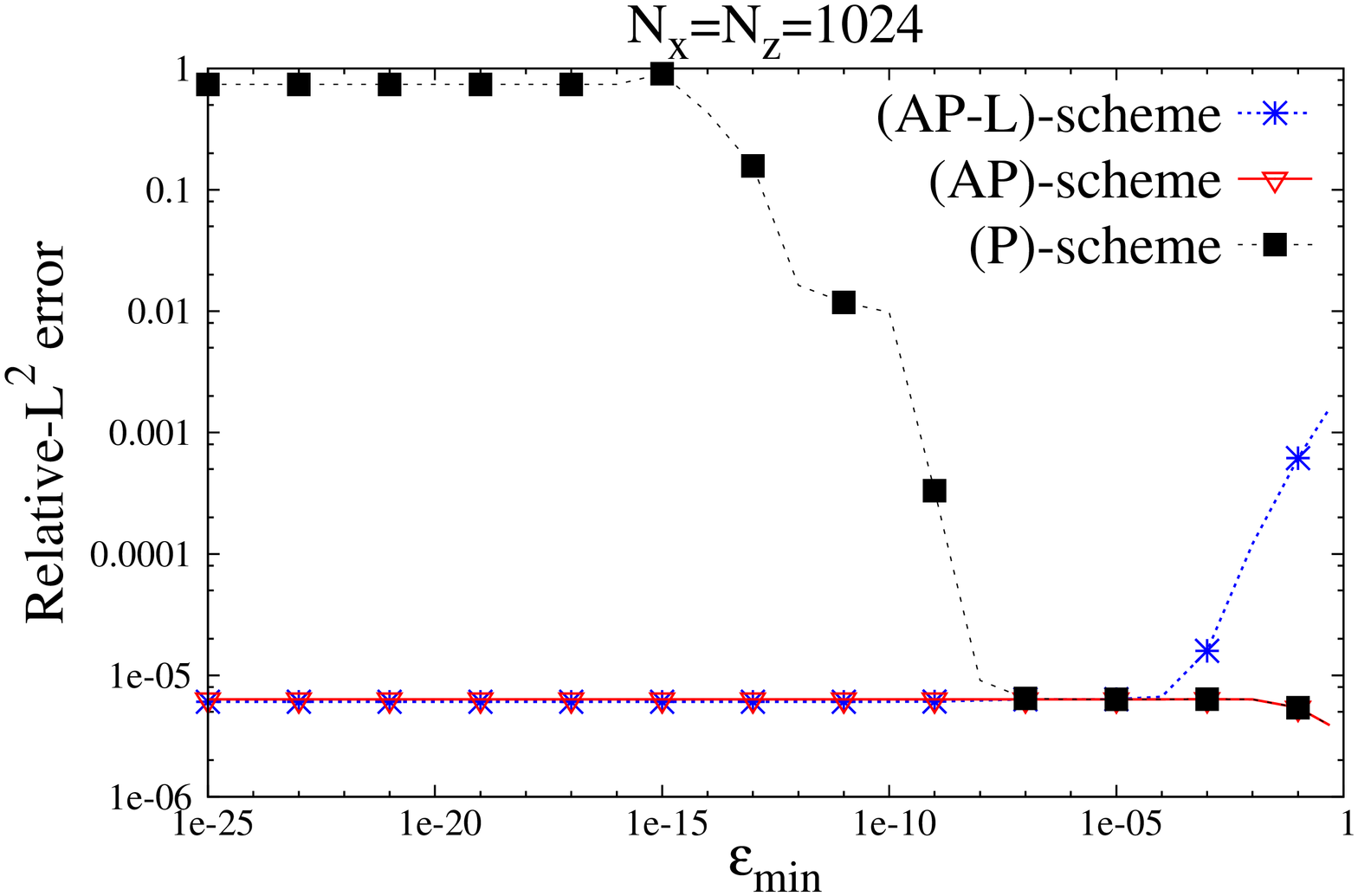}}
\subfigure[\label{fig:ap1024:b}]{\includegraphics[angle=-90,width=0.49\textwidth]{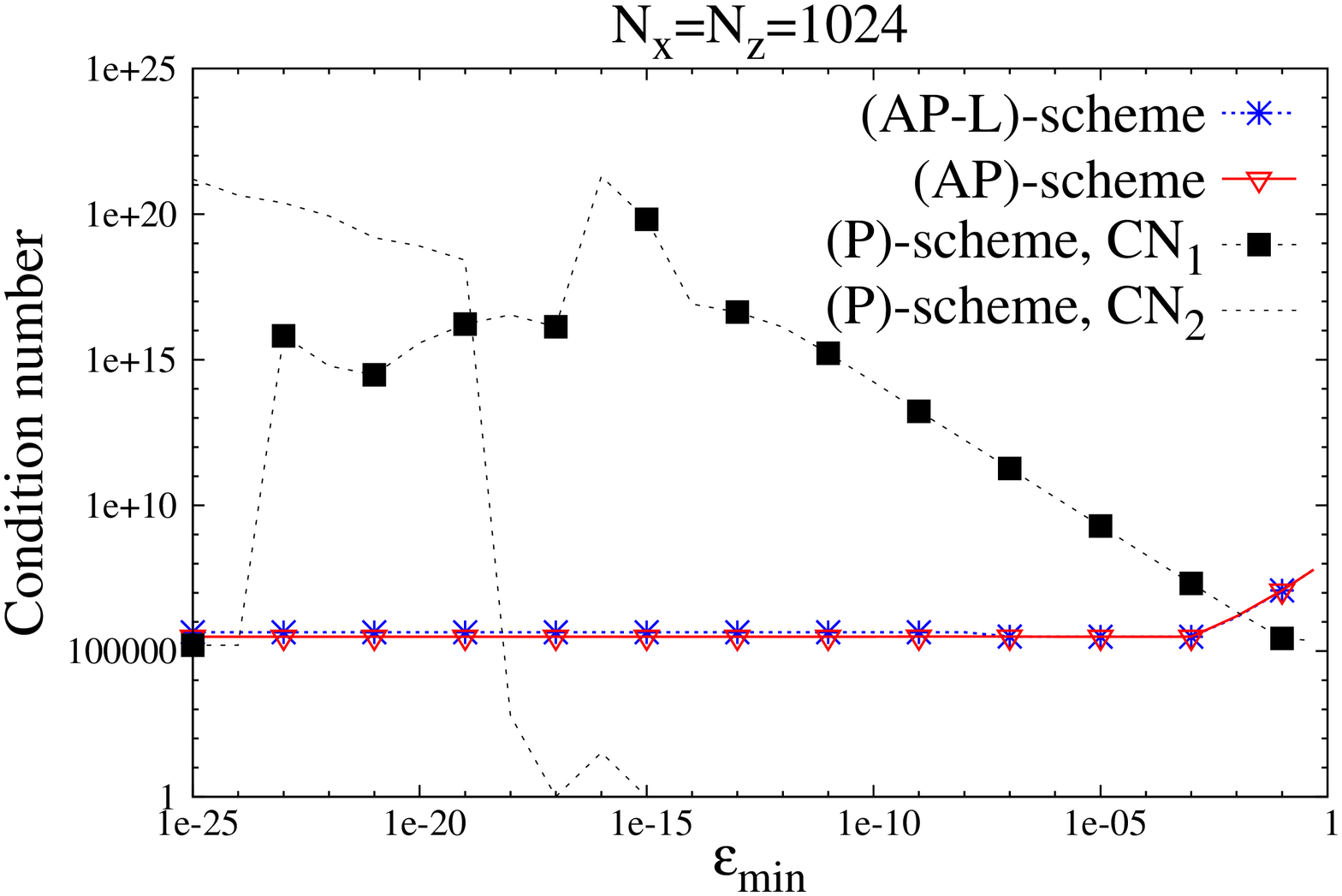}}
\caption{Comparison of the hybrid (AP/L)-model, the AP-scheme and the singular perturbation P-model accuracy for heterogeneous anisotropy ratios as defined by \eqref{eq:defeps} on a $1024\times1024$ mesh: (a) relative L$^2$-error between the exact solution and its numerical approximations and (b) condition number of the linear systems (estimated by MUMPS) as a functions of $\varepsilon_{min}$ (with $\varepsilon_{max}=1$).}
\label{fig:ap1024}
\end{center}
\end{figure}

The condition number of the singular perturbation problem increases almost linearly with vanishing $\varepsilon_\text{min}$. When the condition number is comparable to the computer arithmetic precision, the accuracy of the solution is significantly altered as depicted on Figures~\ref{fig:ap64:a} and \ref{fig:ap1024:a}. On the coarsest mesh the accuracy of the computations is altered for $\varepsilon_\text{min}$-values smaller than $10^{-11}$ (see Figure~\ref{fig:ap64:b}). For these values the condition number of the linear system is reported by the solver to be as large as $10^{13}$. For the computations carried out on the refined mesh, this threshold is reached for larger values of the $\varepsilon_\text{min}$, precisely $\varepsilon_\text{min}=10^{-7}$, corresponding to an estimate of condition number value above $10^{12}$. On this refined mesh, the solution cannot be accurately approximated via the singular perturbation problem for variations of the anisotropy ratio larger than $10^{7}$, as depicted on Figure~\ref{fig:ap1024:a}. The range of anisotropy variations tractable by the P-model is getting narrower with the mesh refinement.

The AP property on the standard AP-scheme translates into a condition number as well as a precision almost independent of the anisotropy strength.

As for the hybrid model, the condition number is observed to be very similar to that of the standard AP-scheme, with almost no variations with respect to the anisotropy strength. Concerning the precision, two regimes can be identified. For small enough  $\varepsilon_\text{min}$-values the accuracy of the computations are comparable to the standard AP-scheme. For computations with $\varepsilon_\text{min}$-values that are not so small (see the growth of the error on Figures~\ref{fig:ap64:a} and \ref{fig:ap1024:a} for the largest values of $\varepsilon_\text{min}$), the limit regime is used in a sub-domain where the asymptotic parameter is too large to guarantee a good approximation of the solution. The approximation error explained by the use of the limit problem in one sub-domain is thus larger than the numerical error of the discretization, producing computations with a poor accuracy. Since the precision of the space discretization increases with vanishing mesh sizes, the accuracy of the computations carried out by the hybrid model is improved by a mesh refinement if the coupling interface is immersed in an area where the value of the asymptotic parameter is small enough. For the coarsest mesh, this requirement is met for $\varepsilon_\text{min}$ values smaller than $10^{-3}$ while this threshold equals $10^{-4}$ for the refined mesh.

%%%%%%%%%%%%%%%%%%%%%%%%%%%%%%%%%%%%%%
\subsubsection{Efficiency of the hybrid numerical method}
%%%%%%%%%%%%%%%%%%%%%%%%%%%%%%%%%%%%%%

Finally, the numerical efficiency of the hybrid method (in terms of simulation time and memory usage) is compared to that of the fully AP-scheme and the discretized singular perturbation problem. Two different sets of computations have been performed. The first one is related to an anisotropy ratio defined by \eqref{eq:defeps} with $\varepsilon_\text{max}=1$ and $\varepsilon_\text{min}=10^{-8}$ and with a domain decomposition verifying $|\Omega^1_z|= \nicefrac{2}{5}\ |\Omega_z|$ ($|\cdot|$ denoting the length) and $|\Omega^2_z|= \nicefrac{3}{5}\ |\Omega_z|$. For the second setup the sub-domain  $\Omega^2_z$ is enlarged with a decomposition yielding to $|\Omega^1_z|= \nicefrac{3}{10}\ |\Omega_z|$ and $|\Omega^2_z|= \nicefrac{7}{10}\ |\Omega_z|$. The anisotropy ratio obeys the same definition \eqref{eq:defeps} but with $\varepsilon_\text{max}=10^{-5}$ and $\varepsilon_\text{min}=10^{-10}$. The main characteristics of the computations are gathered in Table~\ref{tab:time-8} for the first setup and in Table~\ref{tab:time-10} for the second one.%
\begin{table}[!ht]
\begin{center}
%\begin{tabular}{|c|c|c|r|r|r|r|}
\resizebox{\textwidth}{!}{\begin{tabular}{|c|c|r|r|r|r|r|r|}
\hline
Scheme %& $\Omega$ 
& $N_x=N_z$ & Time (AMF) & \#entries in factor (AMF) & \#rows & \#non zeros & L$^2$-error \\
\hline
\hline
AP/L %& $\Omega_b$ 
& 250 & 83\%   (72\%) & 3 705 458 (2 859 836) & 26 000 & 533 324 & $1.06\times 10^{-4}$ \\
\hline
AP %& $\Omega_b$ 
& 250 & 247\%   (187\%) & 10 742 572 (9 301 072) & 63 500 & 1 318 724 & $1.06\times 10^{-4}$ \\
\hline
P %& $\Omega_b$ 
& 250 & 100\%  (100\%) & 5 614 870 (5 726 236) & 63 000 & 563 992 & $1.06\times 10^{-4}$ \\
\hline
\hline
AP/L %& $\Omega_b$ 
& 500 & 83\%   (72\%) & 18 461 132 (13 736 886) & 102 000 & 2 116 674 & $2.64\times 10^{-5}$ \\
\hline
AP %& $\Omega_b$ 
& 500 & 311\%  (219\%) & 61 429 052 (40 562 916) & 252 000 & 5 262 474 & $2.65\times 10^{-5}$ \\
\hline
P %& $\Omega_b$ 
& 500 & 100\%  (100\%) & 26 697 996 (26 940 422) & 251 000 & 2 252 992 & $2.70\times 10^{-5}$ \\
\hline
\hline
AP/L %& $\Omega_b$ 
& 1000 & 72\%   (61\%) & 85 235 506 (61 449 488) & 404 000 & 8 433 374 & $6.52\times 10^{-6}$ \\
\hline
AP %& $\Omega_b$ 
& 1000 & 357\%  (213\%) & 319 347 622 (183 083 450) & 1 004 000 & 21 024 974 & $6.64\times 10^{-6}$ \\
\hline
P %& $\Omega_b$ 
& 1000 & 100\% (100\%) & 130 031 284 (121 110 920) & 1 002 000 & 9 005 992 & $2.23\times 10^{-5}$ \\
\hline
\hline
AP/L %& $\Omega_b$ 
& 2000 & 74\%  (43\%) & 527 455 320 (289 307 676) & 1 608 000 & 33 666 774 & $1.54\times 10^{-6}$ \\
\hline
AP %& $\Omega_b$ 
& 2000 & 351\% (146\%) & 1 678 279 706 (802 178 884) & 4 008 000 & 84 049 974 & $1.66\times 10^{-6}$ \\
\hline
P %& $\Omega_b$ 
& 2000 & 100\% (100\%) & 600 370 134 (557 859 738) & 4 004 000 & 36 011 992 & $8.88\times 10^{-5}$ \\
%%%%%
\hline
\end{tabular}}
\caption{Efficiency of the Hybrid and AP methods compared to the discretized singular perturbation model: computational time relative to that of the P-model with the MUMPS solver and the METIS (AMF) ordering, number of non zero elements after factorization with the METIS (AMF) ordering, number of rows and of non zero elements in the matrices, and precision of the computations  (relative error L$^2$-norm) carried out with  $|\Omega^1_z|= \nicefrac{2}{5} \ |\Omega_z|$ and $|\Omega^2_z|= \nicefrac{3}{5}\ |\Omega_z|$ on different mesh resolutions.}
\label{tab:time-8}
\end{center}
\end{table}%
\begin{table}[!ht]
\begin{center}
%\begin{tabular}{|c|c|c|r|r|r|r|}
\resizebox{\textwidth}{!}{\begin{tabular}{|c|c|r|r|r|r|r|}
\hline
Scheme %& $\Omega$ 
& $N_x=N_z$ & Time (AMF) & \#entries in factor (AMF) & \#rows & \#non zeros & L$^2$-error \\
\hline
\hline
AP/L %& $\Omega_b$ 
& 250 & 57\% (56\%) & 2 431 654 (1 999 636) & 19 750 & 402 424 & $5.61\times 10^{-5}$ \\
\hline
AP %& $\Omega_b$ 
& 250 & 232\% (213\%) & 10 737 412 (9 301 072) & 63 500 & 1 318 724 & $5.61\times 10^{-5}$ \\
\hline
P %& $\Omega_b$ 
& 250 & 100\% (100\%) & 5 614 870 (5 728 812) & 63 000 & 563 992 & $1.51\times 10^{-4}$ \\
\hline
\hline
AP/L %& $\Omega_b$ 
& 500 & 58\% (53\%) & 13 975 832 (9 861 960) & 77 000 & 1 592 374 & $1.41\times 10^{-5}$ \\
\hline
AP %& $\Omega_b$ 
& 500 & 301\% (203\%) & 61 469 878 (40 655 248) & 252 000 & 5 262 474 & $1.41\times 10^{-5}$ \\
\hline
P %& $\Omega_b$ 
& 500 & 100\% (100\%) & 26 699 956 (26 940 422) & 251 000 & 2 252 992 & $3.35\times 10^{-3}$ \\
\hline
\hline
AP/L %& $\Omega_b$ 
& 1000 & 43\% (39\%) & 60 103 598 (45 148 064) & 304 000 & 6 334 774 & $3.58\times 10^{-6}$ \\
\hline
AP %& $\Omega_b$ 
& 1000 & 329\% (190\%) & 320 132 040 (183 805 436) & 1 004 000 & 21 024 974 & $3.53\times 10^{-6}$ \\
\hline
P %& $\Omega_b$ 
& 1000 & 100\% (100\%) & 130 031 284 (121 110 920) & 1 002 000 & 9 005 992 & $1.21\times 10^{-2}$ \\
\hline
\hline
AP/L %& $\Omega_b$ 
& 2000 & 41\% (26\%) & 351 816 570 (206 531 976) & 1 208 000 & 25 269 574 & $1.17\times 10^{-6}$ \\
\hline
AP %& $\Omega_b$ 
& 2000 & 299\% (137\%) & 1 682 372 230 (804 867 106) & 4 008 000 & 84 049 974 & $8.84\times 10^{-7}$ \\
\hline
P %& $\Omega_b$ 
& 2000 & 100\% (100\%) & 600 370 134 (557 859 738) & 4 004 000 & 36 011 992 & $1.50\times 10^{-2}$ \\
%%%%%
\hline
\end{tabular}}
\caption{Efficiency of the Hybrid and AP methods compared to the discretized singular perturbation model: computational time relative to that of the P-model with the MUMPS solver and the METIS (AMF) ordering, number of non zero elements after factorization with the METIS (AMF) ordering, number of rows and of non zero elements in the matrices, and precision of the computations  (relative error L$^2$-norm) carried out with  $|\Omega^1_z|= \nicefrac{3}{10} \ |\Omega_z|$ and $|\Omega^2_z|= \nicefrac{7}{10}\ |\Omega_z|$ on different mesh resolutions.}
\label{tab:time-10}
\end{center}
\end{table}
 Columns 5 and 6 display the number of rows and the number of non zero elements in the matrices stemming from the discretization of the three methods. The computational time required to solve these linear systems using the same solver (MUMPS) with two different matrix ordering algorithms (METIS \cite{METIS} and AMF \cite{AMF}) is collected in column 3. The entries of this column are relative to the computational time of the P-model. The number of non zero elements in the factorized matrix is precised in column 4 for the two matrix orderings mentioned above. The last column gives the relative L$^2$-error between the exact solution and the numerical approximations. The numerical method is specified in the first column and the mesh size in the second one.

The AP-scheme requires a larger memory amount than the P-model to store the linear system matrix. This growth is marginally explained by the increase in the number of unknowns due to the computation of the solution mean part, on the one hand, and to the introduction of the Lagrangian, on the other hand. Fortunately, both unknowns only depend on the $x$ coordinate and explain thus a moderate increase in the number of the matrix rows~: from $N_x\times (N_z+2)$ for the P-model to ($N_x\times (N_z+4)$) for the AP-scheme. However, the AP-scheme matrix exhibits an amount of non zero elements ($(7N_z+13)(3N_x-2)$, see also \cite[Figure~1]{BDNY}) significantly increased, by a factor of about $\nicefrac{7}{3}\approx 2.33$, compared to that of the P-model ($(3N_z+4)(3N_x-2)$). This is equally explained by the $B_l$ and $B_c$ sub-matrices, introduced in \eqref{eq:def:Matrix:AP} for the discretization of the fluctuation zero mean value constraint, and by the coupling between the mean and the fluctuating parts of the solution by means of the $C_a$ and $C_f$ sub-matrices.
These sub-matrices have a large amount of non zero elements compared to the $A_{xa}$ sub-matrix and explain thus most of the increase in the AP-scheme memory usage. A similar conclusion can be drawn from the comparisons of the  factorized matrices. The fill-in of the factorized AP-matrix ranges from 140\% to 280\% that of the P-model (column~4 of Tables~\ref{tab:time-8} and \ref{tab:time-10}).

The improvement in the hybrid method efficiency is obtained thanks to the reduction of the fluctuation computational domain size. This offers a significant decrease in the row number of the hybrid (AP/L)-model matrix \eqref{eq:system:hybrid} as well as a reduction in the size of the sub-matrices with the poorest sparsity. Compared to the AP-scheme, the number of non zero elements stored in the (AP/L)-matrix is roughly divided by 2.4 for the first setup (Table~\ref{tab:time-8}) and by 3.3 for the second one (Table~\ref{tab:time-10}). Compared to the P-model, a gain in the matrix size of the (AP/L)-scheme can be anticipated as soon as the interface position verifies $\iota > \nicefrac{4}{7}N_z\approx0.57N_z$. The fill-in is reduced by 5\% to 9\% for the first setup ($\iota=\nicefrac{3}{5}N_z$) and 20\% to 30\% for the second one ($\iota=\nicefrac{7}{10}N_z$). Importantly also for the linear system solver used for these investigations, this advantage is amplified for the sparsity of the factorized matrix. For the most favorable setup, the fill-in of the (AP/L) factorized matrix is only 43\% to 59\% that of the P-model matrix for the METIS ordering and 26\% to 56\% for the AMF algorithm. This gain results in a significant reduction of the computation time with a speed-up, relatively to the P-model, ranging from $1.2$ to $3.8$ over the data collected in these investigations.

Finally the last column of Tables~\ref{tab:time-8} and \ref{tab:time-10} illustrates the precision of the different methods. For the first setup, both AP-methods produce accurate computations with a small advantage for the hybrid method (at most, less than 10\%, for the largest system size). %This is explained by the smaller size of the hybrid matrix which improves slightly the system conditioning number. 
For the second setup, the hybrid (AP/L)-method fails to reach the accuracy of the standard AP-method for the most refined meshes. The mesh refinement increases the precision of the space discretization and necessitates, in order to improve the overall method accuracy, an interface immersed in a region with small enough $\varepsilon$-values. This requirement is met for the first setup, but not for the two most refined meshes with the second one. Note that the accuracy of the hybrid (AP/L)-method may be recovered by moving the interface location into a region with smaller $\varepsilon$-values.

% The weakness of the P-model are also emphasized by these numerical investigations. In the first setup, the P-model accuracy is observed to deteriorate with the increase of the mesh size. The precision of the computations carried out by the P-model are comparable to that of the (AP) methods only for the two coarsest meshes. This is another illustration of the bad condition number of the P-model matrix. The plots of this number, estimated with a different setup, on Figures~\ref{fig:ap64:b} and \ref{fig:ap1024:b} show a growth of its value with the mesh size, with for instance, a value of about $10^7$ for $\varepsilon_\text{min}=10^{-5}$ on a $64\times64$ mesh and more than $10^9$ on the $1024\times1024$ mesh. A similar deterioration explains the poor accuracy of the P-model reported in Table~\ref{tab:time-8} for the most refined meshes. For the second setup, the anisotropy is more severe in the whole domain making the computations not tractable at all.

\bigskip

These numerical investigations demonstrate a significant advantage of the hybrid method as soon as the limit model can be used on a sub-domain large enough. For the numerical investigations performed in this section, the overhead of the AP-scheme is completely erased. In favorable situations, the hybrid method is more efficient in terms of storage requirements as well as in the computational time than the standard P-model. Of course, these results should be worsen or improved accordingly to the domain decomposition authorized by the anisotropy ratio values and the precision required.

%%%%%%%%%%%%%%%%%%%%%%%%%%%%%%%%%%%%%%%%%%%%%%%
\section{Conclusion}\label{SEC5} 
%%%%%%%%%%%%%%%%%%%%%%%%%%%%%%%%%%%%%%%%%%%%%%%

In this paper, we studied a 2D highly anisotropic elliptic problem, whose initial formulation (called P-model) leads to an ill-conditioned system when the anisotropy parameter $\varepsilon$ tends to zero. Degond \textit{et al.} developed recently an asymptotic-preserving AP-model to overcome this problem (see \cite{DDN} for example). It is uniformly accurate for any value of $\varepsilon\in\left[0,1\right]$, but more costly than the P-model. 

In the context of ionospheric plasmas, $\varepsilon$ appears to be very small in a large region of the domain, such that, in this region, the use of the 1D limit L-model can be interesting from a computational point of view. The coupling of the L-model to the P-model is not always possible, because their validity domains may not overlap.
That is the reason why we proposed to solve the L-model in the region where $\varepsilon$ is close to zero, while an AP-model is used in the rest of the domain. The two models are coupled thanks to Dirichlet-Neumann transfer conditions.%Dirichlet-to-Neumann interface conditions.

The coupling strategy has been detailed and mathematically analyzed. Then, we have presented the numerical discretization and applied our new method in some test cases.
The obtained (AP/L)-model is shown to be accurate and well-conditioned even if $\varepsilon$ presents steep gradients with values closed to zero in one part of the domain and of order one in the other. Moreover, the computational time is significantly reduced compared to the full AP-model, due to the use of the 1D L-model in a large part of the domain.

~

\paragraph{Acknowledgement.} F. Deluzet and C. Negulescu would like to thank Pierre Degond for having initiated this series of works on highly anisotropic elliptic equations. 
This work has been supported by the Agence Nationale de la Recherche (ANR) under contract IODISSEE (IOnospheric DIsturbanceS and SatEllite-to-Earth communications, ANR-09-COSI-007-02) and the FR-FCM (F\'ed\'eration de Recherche pour la Fusion par Confinement Magn\'etique).

%%%%%%%%%%%%%%%%%%%%%%%%%%%%%%%%%%%%%%%%%%%%%%%

%%%%%%%%%%%%%%
\end{document}